# SPECIFICATION TESTING IN NONLINEAR AND NONSTATIONARY TIME SERIES AUTOREGRESSION[1]

By Jiti Gao, Maxwell King, Zudi Lu and Dag Tjøstheim[2]

*University of Adelaide, Monash University, Curtin University of Technology and The University of Adelaide and University of Bergen*

This paper considers a class of nonparametric autoregressive models with nonstationarity. We propose a nonparametric kernel test for the conditional mean and then establish an asymptotic distribution of the proposed test. Both the setting and the results differ from earlier work on nonparametric autoregression with stationarity. In addition, we develop a new bootstrap simulation scheme for the selection of a suitable bandwidth parameter involved in the kernel test as well as the choice of a simulated critical value. The finite-sample performance of the proposed test is assessed using one simulated example and one real data example.

**1. Introduction.** Time series regression analysis has a long history. There have been many studies in using parametric linear autoregressive moving average models [Brockwell and Davis (1990)], parametric nonlinear time series models [see, e.g., Tong (1990), Granger and Teräsvirta (1993)], and nonparametric and semiparametric time series models [Tong (1990), Fan and Yao (2003) and Gao (2007)]. In many existing studies, particularly in the nonparametric situation, the focus of attention has been on the case where the observed time series satisfies a type of stationarity. Such a stationarity assumption is quite restrictive in many cases.

In the parametric time series case, estimation and specification testing methods have been developed to deal with nonstationarity. In recent years, attempts have also been devoted to the estimation of nonlinear and nonstationary time series models using nonparametric methods. Existing studies

Received July 2008; revised February 2009.

[1]Supported by the Australian Research Council Discovery Grants DP-05-58602, DP-08-79088 and DP-09-84686.

[2]Supported by the Danish Research Council Grant 2114-04-0001.

include Phillips and Park (1998) and Karlsen and Tjøstheim (1998, 2001) on nonparametric autoregression, Park and Phillips (2001) on parametric nonlinear regression, Bandi and Phillips (2003) on nonparametric estimation of nonstationary diffusion models, Wang and Phillips (2009) on nonparametric kernel estimation of random walk processes, and Karlsen, Myklebust and Tjøstheim (2007) on nonparametric cointegration. In the original version of this paper, Gao et al. (2006) discuss specification testing problems for both autoregression and conintegration cases with nonstationarity.

In the field of model specification testing with nonstationarity, there is a huge literature on various unit root tests in the parametric linear autoregressive case. To the best of our knowledge, there seems to be very little work on specification testing in the nonparametric nonlinear autoregressive case. This paper aims to discuss such issues. Consider a class of nonlinear autoregressive models of the form

$$X_t = g(X_{t-1}) + u_t, \qquad t = 1, 2, \dots, T, \tag{1.1}$$

where $g(\cdot)$ is an unknown function defined over $R^1 = (-\infty, \infty)$, $\{u_t\}$ is a sequence of independent and identically distributed i.i.d. errors with mean zero and finite variance $\sigma_u^2 = E[u_1^2]$, and $T$ is the number of observations. The initial value $X_0$ of $X_t$ may be any $O_p(1)$ random variable. However, we set $X_0 = 0$ to avoid some unnecessary complications in exposition.

When $g(X_{t-1}) = X_{t-1} + g_1(X_{t-1})$ with $g_1(\cdot)$ being an identifiable nonlinear function, model (1.1) becomes a nonlinear random walk model. Granger, Inoue and Morin (1997) discuss some parametric cases for this model, and suggest several estimation procedures. As $g(\cdot)$ usually represents some kind of nonlinear fluctuation in the conditional mean, it would be both theoretically and practically useful to test whether such a nonlinear term is significant before using model (1.1) in practice. We therefore propose testing the following null hypothesis:

$$H_0 : P(g(X_{t-1}) = X_{t-1}) = 1 \qquad \text{for all } t \geq 1. \tag{1.2}$$

The main difference between our approach and existing ones is that we need not prespecify $g(x)$ parametrically as $g(x) = \theta x$ and then test $H_0' : \theta = 1$ as has been done in the literature. Our approach is that we test $H_0$ nonparametrically. In doing so, we can avoid possibly misspecifying the true model before using a specification testing procedure.

The main contributions of this paper are as follows:

(i) It proposes a nonparametric kernel test for nonlinear nonstationarity against nonlinear stationarity in model (1.1). This test procedure corresponds to the well-known test proposed by Dickey and Fuller (1979) for the parametric case.

(ii) It establishes an asymptotically normal test for testing the conditional mean in model (1.1) under the null hypothesis. Theoretical properties for the proposed test procedure are established.

(iii) This paper is then concerned with discussing the power function of the proposed test under a stationary alternative. Some asymptotic consistency results under both the null and alternative hypotheses are established.

(iv) In order to implement the proposed test in practice, we develop a new simulation procedure based on the assessment of both the size and power functions of the proposed test.

The rest of the paper is organized as follows. Section 2 establishes a simple nonparametric test and an asymptotic distribution under the null hypothesis. Discussion about the power function of the proposed test is given in Section 3. Section 4 shows how to implement the proposed test in practice. Section 5 concludes the paper with some remarks on extensions. Mathematical details are relegated to the Appendix. Some additional derivations are given in Appendices B–E of Gao et al. (2008).

**2. Nonparametric unit root test.** Consider model (1.1) and a general testing problem of the form

$$\begin{aligned} &H_0 : P(g(X_{t-1}) = X_{t-1}) = 1 \quad \text{against} \\ &H_1 : P(g(X_{t-1}) = X_{t-1} + \Delta_T(X_{t-1})) = 1, \end{aligned} \tag{2.1}$$

where $\{\Delta_T(x)\}$ is a sequence of unknown functions.

Before proposing our test statistic for (2.1), we consider the conventional Nadaraya–Watson (NW) kernel estimate of the form

$$\widehat{g}(x) = \sum_{s=1}^{T} W_T(x, X_{s-1}) X_s = \frac{\sum_{s=1}^{T} K_h(X_{s-1} - x) X_s}{\sum_{t=1}^{T} K_h(X_{t-1} - x)}, \tag{2.2}$$

where $W_T(x, X_{s-1}) = \frac{K_h(X_{s-1}-x)}{\sum_{t=1}^{T} K_h(X_{t-1}-x)}$, in which $K_h(\cdot) = K(\cdot/h)$, $K(\cdot)$ is a probability kernel function and $h$ is a bandwidth parameter.

Let $A(X_{t-1}, X_{s-1}) = \frac{1}{T}\sum_{k=1}^{T} W_T(X_{k-1}, X_{t-1}) W_T(X_{k-1}, X_{s-1})$ and $\widehat{X}_{t-1} = \sum_{s=1}^{T} W_T(X_{t-1}, X_{s-1}) X_{s-1}$. We then compare $\widehat{g}(X_{t-1})$ with $\widehat{X}_{t-1}$ by

$$\begin{aligned} N_T(h) = N_T(X_1, \dots, X_T; h) &= \frac{1}{T}\sum_{t=1}^{T} [\widehat{g}(X_{t-1}) - \widehat{X}_{t-1}]^2 \\ &= \sum_{s=1}^{T}\sum_{t=1}^{T} \left( \frac{1}{T} \sum_{k=1}^{T} W_T(X_{k-1}, X_{t-1}) W_T(X_{k-1}, X_{s-1}) \right) u_t u_s \\ &= \sum_{s=1}^{T}\sum_{t=1}^{T} A(X_{t-1}, X_{s-1}) u_t u_s, \end{aligned}$$

where $u_t = X_t - X_{t-1}$ under $H_0$. Similar forms have been used for the stationary time series case [see, e.g., Hjellvik, Yao and Tjøstheim (1998)]. Other alternatives to $N_T(h)$, including the introduction of $M_T(h)$ below, are discussed in Gao et al. (2006).

In theory, we can derive a test statistic based on $N_T(h)$. As can be seen, $N_T(h)$ involves both a triple summation and a kind of random denominator problem, which may cause more difficulty and technicality than those for the stationary case. Compared with $M_T(h)$ below, our experience with the stationary case also shows that a test statistic based on $N_T(h)$ has less attractive properties than those based on $M_T(h)$ below [see, e.g., Li (1999), Gao and King (2004) and Chapter 3 of Gao (2007)].

We thus propose using a test statistic of the form

$$M_T = M_T(h) = \sum_{t=1}^{T} \sum_{s=1, s\neq t}^{T} \widehat{u}_s K_h(X_{s-1} - X_{t-1}) \widehat{u}_t, \tag{2.3}$$

where $\widehat{u}_t = X_t - \widehat{g}(X_{t-1})$. We now introduce the following conditions.

Assumption 2.1. (i) Suppose that $\{u_t\}$ is a sequence of independent and identically distributed i.i.d. errors with $E[u_1] = 0$ and $E[u_1^2] = \sigma_u^2 < \infty$. Let $0 < \mu_4 = E[u_1^4] < \infty$.

(ii) Suppose that $\{u_t\}$ has a symmetric density function $f(u)$. Let $f'(u)$ be the first derivative of $f(u)$ and $f'(u)$ be continuous at $u \in (-\infty, \infty)$. Let $\psi(\cdot)$ be the characteristic function of $\{u_t\}$ satisfying $\int_{-\infty}^{\infty} |v||\psi(v)|\, dv < \infty$.

(iii) Let $K(\cdot)$ be a symmetric probability density function. Suppose that there are constants $c_1 > 0$ and $0 < c_2 < c_3 < \infty$ such that $c_2 I(|u| \le c_1) \le K(u) \le c_3 I(|u| \le c_1)$. In addition, suppose that $|K(x+y) - K(x)| \le \Psi(x)|y|$ for all $x \in C(K)$ and any small $y$, where $\Psi(x)$ is nonnegative bounded function for all $x \in C(K)$ and $C(K)$ denotes the compact support of $K(\cdot)$.

(iv) Assume that $h$ satisfies $\lim_{T\to\infty} T^{3/10} h = 0$ and $\limsup_{T\to\infty} T^{1/2-\varepsilon_0} h = \infty$ for all $0 < \varepsilon_0 < \frac{1}{5}$.

Remark 2.1. The i.i.d. assumption in Assumption 2.1(i) is needed to ensure that the partial sum $S_t = \sum_{s=1}^{t} u_s$ has independent increments, although, $\{S_t\}$ itself is nonstationary and dependent. Under this assumption, we are able to establish the main results of this paper in Theorems 2.1, 2.2 and 3.1 below. Assumption 2.1(ii) imposes some mild conditions on both the density function and the characteristic function and it holds in many cases. The condition $\int_{-\infty}^{\infty} |v||\psi(v)|\, dv < \infty$ is to ensure certain convergence results. Let $\phi_T(x)$ be the density function of $\frac{1}{\sqrt{T}\sigma_u} \sum_{t=1}^{T} u_t$. Then under Assumption 2.1(ii),

$$\sup_x |\phi_T(x) - \phi(x)| \to 0 \quad \text{and} \quad \sup_x |\phi_T'(x) - \phi'(x)| \to 0, \tag{2.4}$$

where $\phi_T'(x)$ and $\phi'(x)$ are first derivatives, and $\phi(x) = \frac{1}{\sqrt{2\pi}} e^{-x^2/2}$ is the density function of the standard normal random variable $N(0,1)$. The proof of (2.4) is quite standard [Chapters 8 and 9 of Chow and Teicher (1988)].

Assumption 2.1(iii) also holds in many cases. For example, when $K(x) = \frac{1}{2} I_{[-1,1]}(x)$, Assumption 2.1(iii) holds automatically. In addition, Assumption 2.1(iv) does not look unnatural in the nonstationary case, although it looks more restrictive than that required for the stationary case. In addition, the conditions of Theorems 5.1 and 5.2 of Karlsen and Tjøstheim (2001) imposed on $h$ become simplified since we are interested in the special case of random walk with a tail index of $\beta = \frac{1}{2}$ involved in those conditions. Such conditions on the bandwidth for nonparametric testing in the nonstationary case are equivalent to the minimal conditions: $\lim_{T\to\infty} h = 0$, $\lim_{T\to\infty} Th = \infty$ and $\lim_{T\to\infty} Th^4 = 0$ required in nonparametric kernel testing for the stationary time series cases [see, e.g., Gao (2007)].

Let $\widehat{\sigma}_T^2 = \widehat{\sigma}_T^2(h) = 2\sum_{t=1}^T \sum_{s=1,s\neq t}^T \widehat{u}_s^2 K_h^2(X_{s-1} - X_{t-1})\widehat{u}_t^2$. As can be seen from the proof of Theorem 2.2 below, under $H_0$ we have for the normalized test statistic

$$\begin{aligned} \widehat{L}_T(h) = \frac{M_T(h)}{\widehat{\sigma}_T(h)} &= \frac{\sum_{t=1}^T \sum_{s=1,s\neq t}^T \widehat{u}_s K((X_{t-1} - X_{s-1})/h)\widehat{u}_t}{\sqrt{2\sum_{t=1}^T \sum_{s=1,s\neq t}^T \widehat{u}_s^2 K^2((X_{t-1} - X_{s-1})/h)\widehat{u}_t^2}} \\ &= \frac{\sum_{t=2}^T \sum_{s=1}^{t-1} u_s K((\sum_{j=s+1}^{t-1} u_j + u_s)/h)u_t}{\sqrt{\sum_{t=2}^T \sum_{s=1}^{t-1} u_s^2 K^2((\sum_{j=s+1}^{t-1} u_j + u_s)/h)u_t^2}} + o_P(1). \end{aligned} \tag{2.5}$$

In comparison with existing forms for the stationary case [e.g., (34) of Arapis and Gao (2006)], establishing an asymptotic distribution for $\widehat{L}_T(h)$ becomes nonstandard mainly due to the fact that $\{X_t\}$ is now nonstationary and $\{u_s\}$ is involved in both the argument of $K(\cdot)$ and in a factor multiplying $K(\cdot)$.

Let

$$\sigma_T^2 = \mathrm{E}\left(\sum_{t=1}^T \sum_{s=1,s\neq t}^T u_s K\left(\frac{X_{t-1} - X_{s-1}}{h}\right) u_t\right)^2.$$

Before we study asymptotic properties of $\widehat{L}_T(h)$, we need to evaluate the asymptotic order of $\sigma_T^2$ in Theorem 2.1 below. The proof is given in Lemma A.1 in Appendix below.

THEOREM 2.1. *Consider model (1.1). Assume that Assumption 2.1 holds. Then under* $H_0$

$$\sigma_T^2 = C_{10} T^{3/2} h(1 + o(1)),$$

*where* $C_{10} = \frac{16\sigma_u^4 J_{02}}{3\sqrt{2\pi}}$, *in which* $\sigma_u^2 = E[u_1^2]$ *and* $J_{02} = \int K^2(x)\,dx$.

Note that $\sigma_T^2$ is proportional to $T^{3/2}h$. When $\{X_t\}$ of model (1.1) is stationary, however, $\sigma_T^2$ is proportional to $T^2 h$ as has been given in the literature [Gao (2007)]. Theorem 2.2 below shows that standard normality can still be the limiting distribution of a test statistic under nonstationarity.

THEOREM 2.2. *Consider model (1.1). Suppose that Assumption 2.1 holds. Then under $H_0$ and as $T \to \infty$*

$$\begin{aligned}\widehat{L}_T(h) = \frac{M_T(h)}{\widehat{\sigma}_T(h)} &= \frac{\sum_{t=1}^T \sum_{s=1,s\neq t}^T \widehat{u}_s K_h(X_{s-1}-X_{t-1})\widehat{u}_t}{\sqrt{2\sum_{t=1}^T \sum_{s=1,s\neq t}^T \widehat{u}_s^2 K_h^2(X_{s-1}-X_{t-1})\widehat{u}_t^2}} \\ &= \frac{\sum_{t=2}^T \sum_{s=1}^{t-1} \widehat{u}_s K_h(X_{s-1}-X_{t-1})\widehat{u}_t}{\sqrt{\sum_{t=2}^T \sum_{s=1}^{t-1} \widehat{u}_s^2 K_h^2(X_{s-1}-X_{t-1})\widehat{u}_t^2}} \to_D N(0,1).\end{aligned} \tag{2.6}$$

PROOF. Observe that under $H_0$

$$\begin{aligned}M_T(h) &= \sum_{t=1}^T \sum_{s=1,s\neq t}^T \widehat{u}_s K_h(X_{s-1}-X_{t-1})\widehat{u}_t \\ &= \sum_{t=1}^T \sum_{s=1,s\neq t}^T u_s K_h(X_{s-1}-X_{t-1})u_t \\ &\quad + \sum_{t=1}^T \sum_{s=1,s\neq t}^T \widehat{\delta}_s K_h(X_{s-1}-X_{t-1})\widehat{\delta}_t \\ &\quad + 2\sum_{t=1}^T \sum_{s=1,s\neq t}^T u_s K_h(X_{s-1}-X_{t-1})\widehat{\delta}_t \\ &\equiv M_{T1} + M_{T2} + M_{T3},\end{aligned} \tag{2.7}$$

$$\begin{aligned}\widehat{\sigma}_T^2 &= 2\sum_{t=1}^T \sum_{s=1,s\neq t}^T \widehat{u}_s^2 K_h^2(X_{s-1}-X_{t-1})\widehat{u}_t^2 \\ &= 2\sum_{t=1}^T \sum_{s=1,s\neq t}^T u_s^2 K_h^2(X_{s-1}-X_{t-1})u_t^2 \\ &\quad + 2\sum_{t=1}^T \sum_{s=1,s\neq t}^T \widehat{\delta}_s^2 K_h^2(X_{s-1}-X_{t-1})\widehat{\delta}_t^2 + \widehat{R}_T,\end{aligned} \tag{2.8}$$

where $\widehat{\delta}_t = X_{t-1} - \widehat{g}(X_{t-1})$ and $\widehat{R}_T$ is the remainder term given by

$$\widehat{R}_T = \widehat{\sigma}_T^2 - 2\sum_{t=1}^{T}\sum_{s=1,s\neq t}^{T} u_s^2 K_h^2(X_{s-1} - X_{t-1})u_t^2 - 2\sum_{t=1}^{T}\sum_{s=1,s\neq t}^{T} \widehat{\delta}_s^2 K_h^2(X_{s-1} - X_{t-1})\widehat{\delta}_t^2.$$

In view of (2.7) and (2.8), to prove Theorem 2.2, it suffices to show that as $T \to \infty$

$$\frac{M_{T1}}{\widetilde{\sigma}_T} \to_D N(0,1), \tag{2.9}$$

$$\frac{M_{Ti}}{\widetilde{\sigma}_T} \to_P 0 \qquad \text{for } i = 2, 3, \tag{2.10}$$

$$\frac{\widehat{\sigma}_T^2 - \widetilde{\sigma}_T^2}{\widetilde{\sigma}_T^2} \to_P 0, \tag{2.11}$$

where $\widetilde{\sigma}_T^2 = 2\sum_{t=1}^{T}\sum_{s=1,s\neq t}^{T} u_s^2 K_h^2(X_{s-1} - X_{t-1})u_t^2$.

The proof of (2.9) is given in Lemma A.3 of the Appendix below. In view of (2.9), to complete the proof of Theorem 2.2, it suffices to prove (2.10) and (2.11). We now give the proof of (2.10) and then an outline of the proof of (2.11).

It follows from (2.9) that

$$\frac{1}{\widetilde{\sigma}_T}\sum_{t=1}^{T}\sum_{s=1}^{T} u_s K\left(\frac{X_{t-1} - X_{s-1}}{h}\right)u_t = O_P(1). \tag{2.12}$$

In order to prove (2.10), we first need to show that

$$\frac{M_{T2}}{\sigma_T} = o_P(1). \tag{2.13}$$

Observe that under $H_0 : X_t = X_{t-1} + u_t$

$$\begin{aligned} \widehat{\delta}_t &= X_{t-1} - \widehat{g}(X_{t-1}) \\ &= X_{t-1} - \sum_{s=1}^{T} W_T(X_{t-1}, X_{s-1})X_s \\ &= X_{t-1} - \sum_{s=1}^{T} W_T(X_{t-1}, X_{s-1})X_{s-1} - \sum_{s=1}^{T} W_T(X_{t-1}, X_{s-1})u_s \\ &= \widetilde{X}_{t-1} - \overline{u}_t, \end{aligned} \tag{2.14}$$

where $\widetilde{X}_{t-1} = X_{t-1} - \sum_{s=1}^T W_T(X_{t-1}, X_{s-1})X_{s-1}$ and $\overline{u}_t = \sum_{s=1}^T W_T(X_{t-1}, X_{s-1})u_s$.

Thus, in order to show (2.13), it suffices to show that

$$\sum_{t=1}^T \sum_{s=1,s\neq t}^T \widetilde{X}_{s-1} K\left(\frac{X_{t-1}-X_{s-1}}{h}\right)\widetilde{X}_{t-1} = o_P(\sigma_T), \tag{2.15}$$

$$\sum_{t=1}^T \sum_{s=1,s\neq t}^T \overline{u}_s K\left(\frac{X_{t-1}-X_{s-1}}{h}\right)\overline{u}_t = o_P(\sigma_T). \tag{2.16}$$

The proof of (2.16) is quite technical and thus relegated to Lemma E.1 in Appendix E of Gao et al. (2008). Meanwhile, Assumption 2.1(iii) and a conventional approach [see, e.g., the proof of Theorem 5.1 of Karlsen and Tjøstheim (2001)] imply that uniformly in $x$,

$$\begin{aligned}\widetilde{g}(x) &= g(x) - \sum_{s=1}^T W_T(x, X_{s-1})g(X_{s-1})\\ &= \frac{\sum_{s=1}^T K((x-X_{s-1})/h)(g(x)-g(X_{s-1}))}{\sum_{l=1}^T K((x-X_{l-1})/h)} \\ &= g'(x)h\int uK(u)\,du(1+o_P(1)) = o_P(h),\end{aligned} \tag{2.17}$$

when $g(\cdot)$ is differentiable and the first derivative, $g'(x)$, is continuous.

Using (2.17) for the case of $g(x)=x$, in order to prove (2.15), it suffices to show that

$$h^2 \sum_{t=1}^T \sum_{s=1,s\neq t}^T K\left(\frac{X_{t-1}-X_{s-1}}{h}\right) = o_P(\sigma_T), \tag{2.18}$$

which follows from

$$\sum_{t=1}^T \sum_{s=1,s\neq t}^T E\left[K\left(\frac{X_{t-1}-X_{s-1}}{h}\right)\right] = O(T^{3/2}h) \tag{2.19}$$

and Assumption 2.1(iv). The verification of (2.19) is similar to but simpler than that of (A.3) below.

Hence, (2.13) and (A.50) in the Appendix below imply

$$\frac{M_{T2}}{\widetilde{\sigma}_T} = \frac{M_{T2}}{\sigma_T}\frac{\sigma_T}{\widetilde{\sigma}_T} = o_P(1). \tag{2.20}$$

This proves (2.10) for $i=2$. Furthermore, the proof of (2.10) for $i=3$ follows from (2.12)–(2.20) and

$$\left|\sum_{t=1}^T \sum_{s=1,s\neq t}^T u_s \sqrt{K\left(\frac{X_{t-1}-X_{s-1}}{h}\right)}\sqrt{K\left(\frac{X_{t-1}-X_{s-1}}{h}\right)}\widehat{\delta}_t\right|^2$$

$$\leq \sum_{t=1}^{T} \sum_{s=1,s\neq t}^{T} u_s K\left(\frac{X_{t-1}-X_{s-1}}{h}\right) u_t \sum_{t=1}^{T} \sum_{s=1,s\neq t}^{T} \widehat{\delta}_s K\left(\frac{X_{t-1}-X_{s-1}}{h}\right) \widehat{\delta}_t$$

$$= O_P(\widetilde{\sigma}_T) \cdot o_P(\widetilde{\sigma}_T) = o_P(\widetilde{\sigma}_T^2),$$

where $\widehat{\delta}_t = X_{t-1} - \widehat{g}(X_{t-1})$.

In view of the definitions of $\widehat{\sigma}_T^2$, $\widetilde{\sigma}_T^2$, (2.17) above, (A.50) in the Appendix and

$$\frac{\widehat{\sigma}_T^2 - \widetilde{\sigma}_T^2}{\widetilde{\sigma}_T^2} = \frac{\widehat{\sigma}_T^2 - \widetilde{\sigma}_T^2}{\sigma_T^2} \cdot \frac{\sigma_T^2}{\widetilde{\sigma}_T^2},$$

in order to prove (2.11), it suffices to show that

$$\sum_{t=1}^{T} \sum_{s=1,s\neq t}^{T} \overline{u}_s^2 K^2\left(\frac{X_{s-1}-X_{t-1}}{h}\right) \overline{u}_t^2 = o_P(\sigma_T^2). \tag{2.21}$$

The verification of (2.21) is similar to that of (2.16). This completes the proof of Theorem 2.2. □

Existing studies of test statistics analogous to $\widehat{L}_T(h)$ for the stationary time series case show that the size function of the test is not well approximated using a normal limit distribution. The main reasons are as follows: (a) the rate of convergence of each $\widehat{L}_T(h)$ to asymptotic normality is quite slow even when $\{u_t\}$ is a sequence of independent and identically distributed errors; and (b) the use of a single bandwidth based on an optimal estimation criterion may not be optimal for testing purposes.

In order to improve the finite sample performance of $\widehat{L}_T(h)$, we propose using a bootstrap simulation method. Such a method is known to work quite well in the stationary case. For each given bandwidth satisfying certain theoretical conditions, instead of using an asymptotic critical value of $l_{0.05} = 1.645$ at the 5% level for example, we use a simulated critical value for computing the size function and then the power function. An optimal bandwidth is chosen such that while the size function is controlled by a significance level, the power function is maximized at the optimal bandwidth. Our finite-sample studies show that there is little size distortion when using such a simulated critical value. These issues are discussed in Section 3 below.

**3. Bootstrap simulation and asymptotic theory.** In order to assess the performance of both the size and power function, we need to discuss how to simulate critical values for the implementation of $\widehat{L}_T(h)$ in each case. We then examine the finite sample performance through using two examples in Section 4 below. Before we look at how to implement $\widehat{L}_T(h)$ in practice, we propose the following simulation scheme.

*Simulation scheme*: the exact $\alpha$-level critical value, $l_\alpha(h)$ $(0<\alpha<1)$, is the $1-\alpha$ quantile of the exact finite-sample distribution of $\widehat{L}_T(h)$. Because there are unknown quantities, such as unknown parameters and functions, we cannot evaluate $l_\alpha(h)$ in practice. We propose choosing an approximate $\alpha$-level critical value, $l_\alpha^*(h)$, by using the following simulation procedure:

- Let $X_0=0$. For each $t=1,2,\ldots,T$, generate $X_t^*=X_{t-1}+\widehat{\sigma}_u\varepsilon_t^*$, where $\widehat{\sigma}_u^2$ is a consistent estimator of $\sigma_u^2=E[u_1^2]$ based on the original sample $(X_1,X_2,\ldots,X_T)$, and $\{\varepsilon_t^*\}$ is constructed using either a parametric bootstrap method or a nonparametric bootstrap method.
- Use the data set $\{X_t^*:t=1,2,\ldots,T\}$ to re-estimate $g(\cdot)$ by $\widehat{g}^*(x)=\sum_{s=1}^T W_T\times(x,X_{s-1})X_s^*$. Let $\widehat{u}_t^*=X_t^*-\widehat{g}^*(X_{t-1})$. Compute the test statistic $\widehat{L}_T^*(h)$ that is the corresponding version of $\widehat{L}_T(h)$ by replacing $\widehat{u}_t$ with $\widehat{u}_t^*$ on the right-hand side of $\widehat{L}_T(h)$.
- Repeat the above steps $M$ times and produce $M$ versions of $\widehat{L}_T^*(h)$ denoted by $\widehat{L}_{Tm}^*(h)$ for $m=1,2,\ldots,M$. Use the $M$ values of $\widehat{L}_{Tm}^*(h)$ to construct their empirical bootstrap distribution function. The bootstrap distribution of $\widehat{L}_T^*(h)$ given $\mathcal{X}_T=\{X_t:1\le t\le T\}$ is defined by $P^*(\widehat{L}_T^*(h)\le x)=P(\widehat{L}_T^*(h)\le x|\mathcal{X}_T)$. Let $l_\alpha^*(h)$ satisfy

$$P^*(\widehat{L}_T^*(h)>l_\alpha^*(h))=\alpha$$

  and then estimate $l_\alpha(h)$ by $l_\alpha^*(h)$.
- Define the size and power functions by

$$\alpha(h)=P(\widehat{L}_T(h)\ge l_\alpha^*(h)|H_0) \quad \text{and} \quad \beta(h)=P(\widehat{L}_T(h)\ge l_\alpha^*(h)|H_1).$$

  Let $\mathcal{H}=\{h:\alpha(h)\le\alpha\}$. Choose an optimal bandwidth $\widehat{h}_{\text{test}}$ such that

$$\widehat{h}_{\text{test}}=\arg\max_{h\in\mathcal{H}}\beta(h).$$

  We then use $l_\alpha^*(h_{\text{test}})$ in the computation of both the size and power values of $\widehat{L}_T(\widehat{h}_{\text{test}})$ for each case.

To study the power function of $\widehat{L}_T(h)$, we specify a sequence of alternatives of the form:

$$H_1: P(g(X_{t-1})=X_{t-1}+\Delta_T(X_{t-1}))=1, \tag{3.1}$$

where $\Delta_T(x)$ is a sequence of nonparametrically unknown functions satisfying certain conditions in Assumption 3.2 below.

Under $H_1$, model (1.1) becomes

$$X_t=g(X_{t-1})+u_t=X_{t-1}+\Delta_T(X_{t-1})+u_t, \tag{3.2}$$

where $\Delta_T(x)$ can be consistently estimated by

$$\widehat{\Delta}_T(x) = \frac{\sum_{t=1}^T K_{\widehat{h}_{\mathrm{cv}}}(X_{t-1}-x)(X_t - X_{t-1})}{\sum_{t=1}^T K_{\widehat{h}_{\mathrm{cv}}}(X_{t-1}-x)} \tag{3.3}$$

with $\widehat{h}_{\mathrm{cv}}$ being chosen by a conventional cross-validation selection method.

To establish Theorem 3.1 below, we need the following conditions.

Assumption 3.1. (i) Assumption 2.1 holds.

(ii) Suppose that $g(x)$ is differentiable in $x \in R^1 = (-\infty, \infty)$ and that the first derivative $g'(x)$ is continuous in $x \in R^1$. In addition, $g(x)$ is chosen such that $\{X_t\}$ of (1.1) under $H_1$ is strictly stationary.

Assumption 3.2. Let $f(x)$ be the marginal density function of $\{X_t\}$ under $H_1$. Suppose that $\{\Delta_T(x)\}$ is either an unknown function of the form $\Delta(x)$ or a sequence of unknown functions satisfying

$$\lim_{T\to\infty} T^{5/4}\sqrt{h}\delta^2(T) = \infty \qquad \text{where } \delta^2(T) = \int \Delta_T^2(x) f^2(x)\,dx. \tag{3.4}$$

Since $g(x)$ is not necessarily identical to $x$ under $H_1$, Assumption 3.1(ii) requires that the main interest of this paper is to test linear nonstationarity against nonlinear stationarity. Some secondary conditions on the form of $g(\cdot)$ such that $\{X_t\}$ is strictly stationary under $H_1$ are available from Masry and Tjøstheim (1995).

Assumption 3.2 basically requires that there is some "distance" between $g(X_{t-1})$ and $X_{t-1}$ when $H_0$ is not true. Obviously, there are many different ways of choosing $\Delta_T(x)$ for $H_1$. For example, we may consider testing nonstationarity against stationarity of the form

$$\begin{aligned} &H_0 : X_t = X_{t-1} + u_t \quad \text{versus} \\ &H_1 : X_t = g(X_{t-1}) + u_t = X_{t-1} + \Delta(X_{t-1}) + u_t, \end{aligned} \tag{3.5}$$

where $\{u_t\}$ is a sequence of i.i.d. errors with $E[u_1] = 0$ and $E[u_1^2] = \sigma_u^2 < \infty$, and $\Delta(\cdot)$ can be either a nonparametric or semiparametric function and is chosen such that $\{X_t\}$ is stationary under $H_1$. In this case, we have

$$\begin{aligned} &\frac{1}{T^2 h}\sum_{t=1}^T\sum_{s=1}^T E\left[(g(X_{s-1}) - X_{s-1})K\left(\frac{X_{s-1}-X_{t-1}}{h}\right)(g(X_{t-1}) - X_{t-1})\right] \\ &\qquad = (1+o(1))\int \Delta^2(x) f^2(x)\,dx > 0, \end{aligned} \tag{3.6}$$

since $\{X_t\}$ under $H_1$ is strictly stationary with $f(\cdot)$ being its marginal density function. This, along with Assumption 2.1(iv), implies that Assumption 3.2 holds when $\Delta_T(x) = \Delta(x)$.

We now state the following results and their proofs are given below.

THEOREM 3.1. (i) *Assume that Assumption 3.1 holds. Then under* $H_0$

$$\lim_{T\to\infty} P(\widehat{L}_T(h) > l_\alpha^*) = \alpha.$$

(ii) *Assume that Assumptions 3.1 and 3.2 hold. Then under* $H_1$

$$\lim_{T\to\infty} P(\widehat{L}_T(h) > l_\alpha^*) = 1.$$

Theorems 3.1(i) implies that each $l_\alpha^*$ is an asymptotically correct $\alpha$-level critical value under $H_0$, while Theorem 3.1(ii) shows that $\widehat{L}_T(h)$ is asymptotically consistent against alternatives of the form (3.1) whenever $\delta(T) \geq CT^{-5/8}h^{-1/4}$ for some finite $C > 0$ in this kind of nonparametric testing of nonstationarity against stationarity.

PROOF OF THEOREM 3.1. Recall $\widehat{g}^*(x) = \sum_{s=1}^T W_T(x, X_{s-1})X_s^*$ and $\widehat{u}_t^* = X_t^* - \widehat{g}^*(X_{t-1})$. Let $\widehat{\delta}_t^* = X_{t-1} - \widehat{g}^*(X_{t-1})$. We now have

$$\begin{aligned}
M_T^*(h) &\equiv \sum_{t=1}^T \sum_{s=1,s\neq t}^T \widehat{u}_s^* K_h(X_{s-1} - X_{t-1})\widehat{u}_t^* \\
&= \sum_{t=1}^T \sum_{s=1,s\neq t}^T \widehat{\sigma}_u \varepsilon_s^* K_h(X_{s-1} - X_{t-1})\widehat{\sigma}_u \varepsilon_t^* \\
&\quad + \sum_{t=1}^T \sum_{s=1,s\neq t}^T \widehat{\delta}_s^* K_h(X_{s-1} - X_{t-1})\widehat{\delta}_t^* \qquad (3.7)\\
&\quad + 2\sum_{t=1}^T \sum_{s=1,s\neq t}^T \widehat{\sigma}_u \varepsilon_s^* K_h(X_{s-1} - X_{t-1})\widehat{\delta}_t^* \\
&\equiv M_{T1}^* + M_{T2}^* + M_{T3}^*.
\end{aligned}$$

Using Assumptions 2.1 and 3.1, in view of the notation of $\widehat{L}_T^*(h)$ introduced in the simulation scheme proposed just above Assumption 3.1 as well as the proof of Theorem 2.2, we can show that as $T \to \infty$

$$P^*(\widehat{L}_T^*(h) \leq x) \to \Phi(x) \qquad \text{for all } x \in (-\infty, \infty) \tag{3.8}$$

holds in probability with respect to the distribution of the original sample $\{X_{t-1} : 1 \leq t \leq T\}$, where $\Phi(\cdot)$ is the distribution function of the standard normal random variable $N(0,1)$. In order to prove (3.8), in view of the fact that $\{\varepsilon_s^*\}$ and $\{X_t\}$ are independent for all $s, t \geq 1$, we can show that the proofs of Lemmas A.1 and A.3–A.6 below all remain true by successive conditioning arguments.

Let $z_\alpha$ be the $1-\alpha$ quantile of $\Phi(\cdot)$ such that $\Phi(z_\alpha) = 1-\alpha$. Then it follows from (3.8) that as $T \to \infty$

$$P^*(\widehat{L}_T^*(h) \geq z_\alpha) \to 1 - \Phi(z_\alpha) = \alpha. \tag{3.9}$$

This, together with the construction that $P^*(\widehat{L}_T^*(h) > l_\alpha^*(h)) = \alpha$, implies that as $T \to \infty$

$$l_\alpha^*(h) - z_\alpha \to_P 0. \tag{3.10}$$

Using the conclusion of Theorem 2.2 and (3.8) again, we have that as $T \to \infty$

$$P^*(\widehat{L}_T^*(h) \leq x) - P(\widehat{L}_T(h) \leq x) \to_P 0 \qquad \text{for all } x \in (-\infty, \infty). \tag{3.11}$$

This, along with the construction that $P^*(\widehat{L}_T^*(h) > l_\alpha^*(h)) = \alpha$ again, implies

$$\lim_{T\to\infty} P(\widehat{L}_T(h) > l_\alpha^*(h)) = \alpha. \tag{3.12}$$

Therefore, the conclusion of Theorem 3.1(i) is proved.

Recall $\widehat{u}_t = X_t - \widehat{g}(X_{t-1})$ and let $\lambda_t = X_{t-1} - g(X_{t-1})$. To prove Theorem 3.1(ii), we need to recall the decomposition of $M_T(h)$ in (2.7). Recalling $\widehat{\delta}_t = X_{t-1} - \widehat{g}(X_{t-1})$ and $\lambda_t = X_{t-1} - g(X_{t-1})$, we have

$$\begin{aligned} \widehat{\delta}_t &= X_{t-1} - \widehat{g}(X_{t-1}) = X_{t-1} - g(X_{t-1}) + g(X_{t-1}) - \widehat{g}(X_{t-1}) \\ &= X_{t-1} - g(X_{t-1}) + g(X_{t-1}) - \sum_{s=1}^T W_T(X_{t-1}, X_{s-1}) g(X_{s-1}) \\ &\quad - \sum_{s=1}^T W_T(X_{t-1}, X_{s-1}) u_s = \lambda_t + \widetilde{g}(X_{t-1}) - \overline{u}_t, \end{aligned}$$

where $\widetilde{g}(X_{t-1}) = g(X_{t-1}) - \sum_{s=1}^T W_T(X_{t-1}, X_{s-1}) g(X_{s-1})$. In view of the proof of Theorem 2.2, (2.15)–(2.17) in particular as well as (3.10), in order to prove Theorem 3.1(ii), it suffices to show that under $H_1$

$$\frac{\sum_{s=1}^T \sum_{t=1}^T \lambda_s K_h(X_{s-1} - X_{t-1}) \lambda_t}{\sigma_T} \to_P \infty. \tag{3.13}$$

Similarly to (3.6), we have under $H_1$

$$\begin{aligned} &\frac{1}{\sigma_T} \sum_{t=1}^T \sum_{s=1}^T E\left[(g(X_{s-1}) - X_{s-1}) K\left(\frac{X_{s-1} - X_{t-1}}{h}\right)(g(X_{t-1}) - X_{t-1})\right] \\ &\quad = \frac{T^2 h}{\sigma_T}(1 + o(1)) \int \Delta_T^2(x) f^2(x)\, dx = C T^{5/4} \sqrt{h} \delta^2(T)(1 + o(1)). \end{aligned} \tag{3.14}$$

The verification of (3.13) follows from (3.14) and Assumption 3.2. This finishes the proof of Theorem 3.1. □

Section 4 below shows how to illustrate Theorem 3.1 through using a simulated example and then a real data application.

**4. Examples of implementation.** This section studies some finite-sample properties of both the size and power functions of the proposed test through using two examples. Example 4.1 assesses the finite-sample performance using simulated data. A real data application is given in Example 4.2. Throughout Examples 4.1 and 4.2 below, we use $K(x) = \frac{1}{2} I_{[-1,1]}(x)$.

EXAMPLE 4.1. Consider a nonlinear time series model of the form

$$X_t = X_{t-1} + \Delta(X_{t-1}) + u_t, \tag{4.1}$$

where $X_0 = 0$, $\{u_t\}$ is a sequence of independent normal random errors with $E[u_1] = 0$ and $E[u_1^2] = \sigma_u^2 < \infty$, and $\Delta(x)$ is chosen as a known parametric function with some unknown parameters in the following data generating process.

We then consider two different cases as follows:

$$\begin{aligned} &H_0 : X_t = X_{t-1} + u_t \quad \text{versus} \\ &H_1 : X_t = X_{t-1} + \beta X_{t-1} + u_t \end{aligned} \tag{4.2}$$

and

$$\begin{aligned} &H_0 : X_t = X_{t-1} + u_t \quad \text{versus} \\ &H_1 : X_t = X_{t-1} + \beta X_{t-1} + \frac{\beta}{1 + |X_{t-1}|^\gamma} + u_t, \end{aligned} \tag{4.3}$$

where $0 < \gamma < \infty$, $-2 < \beta < 0$ and $0 < \sigma_u < \infty$ are unknown parameters to be estimated using the conventional MLE method [see Granger and Teräsvirta (1993)].

Since we are interested in assessing the performance of the proposed test for a number of different values for $\beta$, the fixed values of $\sigma_u^2 = 0.05$ and $\gamma = \frac{1}{2}$ are used in generating the data. In addition to the case of $\sigma_u^2 = 0.05$, we have also tried some other values of $\sigma_u$. As our preliminary results show that the resulting finite sample results are very similar, we focus on the case of $\sigma_u^2 = 0.05$ in this example.

Note that $\{X_t\}$ of (4.2) is nonstationary under $H_0$, while it strictly stationary and $\alpha$-mixing under $H_1$ with $0 < \gamma < \infty$, and $-2 < \beta < 0$ in both cases. With the choice of the values for $\beta$ and $\gamma$, the time series $\{X_t\}$ of (4.3) is also strictly stationary under $H_1$ [see, e.g., Masry and Tjøstheim (1995)]. In the simulation, we consider various values of $-2 < \beta < 0$ when computing the power of $\widehat{L}_T(h)$.

As pointed out in the literature for the i.i.d. and stationary time series cases [Hjellvik, Yao and Tjøstheim (1998), Li and Wang (1998), Fan and Linton (2003), Gao (2007) and Gao and Gijbels (2008)], the choice of a kernel bandwidth for testing purposes is quite critical and difficult. In the nonstationary case, however, how to choose an optimal bandwidth parameter is still an open problem.

Thus, in the finite-sample study, we apply the first part of the simulation scheme proposed in Section 3 to simulate a bootstrap critical function $l_\alpha^*(h)$ for each given $h$ in each individual case. We then choose an optimal value for $h$ in each case such that the power function is maximized at such an optimal $\widehat{h}_{\mathrm{test}}$. For each case of $T = 250$, 500 or 750, the finite-sample assessment of the corresponding size and power functions suggests choosing $\widehat{h}_{\mathrm{test}} = 0.160$ when $T = 250$, 0.117 for $T = 500$ and 0.097 when $T = 750$.

To assess the variability of both the size and power with respect to various bandwidth values, we then consider a set of bandwidth values of the form

$$h_i = \frac{1}{2^{5-i}} \widehat{h}_{\mathrm{test}}$$

for $1 \le i \le 5$ with $L_5 = \widehat{L}_T(\widehat{h}_{\mathrm{test}})$. To simplify the notation, we introduce $L_i = \widehat{L}_T(h_i)$ for $1 \le i \le 5$. Since the alternative of model (4.2) is a linear form, we may compare our test with a version of the Dickey–Fuller test of the form [Dickey and Fuller (1979)]

$$L_0 = \frac{\sum_{t=2}^T (X_t - X_{t-1}) X_{t-1}}{\widehat{\sigma}_T \sqrt{\sum_{t=2}^T X_{t-1}^2}}, \tag{4.4}$$

where $\widehat{\sigma}_T^2 = \frac{1}{T} \sum_{t=1}^T (X_t - X_{t-1} - \widehat{\beta}_T X_{t-1})^2$ with $\widehat{\beta}_T = \frac{\sum_{t=2}^T (X_t - X_{t-1}) X_{t-1}}{\sum_{t=2}^T X_{t-1}^2}$.

In the following tables, we consider cases where the number of replications of each of the sample versions of the size and power functions was $M = 1000$, each with $B = 250$ number of bootstrapping resamples $\{\varepsilon_t^*\}$ (involved in the simulation scheme in Section 3 above) from the standard normal distribution $N(0,1)$, and the simulations were done for the cases of $T = 250$, 500 and 750.

Table 1 shows that while the sizes are comparable, the conventional test $L_0$ is more powerful than the proposed test $L_5$ as expected when the alternative model is a linear autoregressive model. However, the biggest power reduction is only about 36% in the case of $T = 250$ and $\beta = -0.05$. This may suggest that we should use the proposed test for nonstationarity in the conditional mean when there is no priori information about the form of the conditional mean.

When the alternative is a nonlinear parametric form as in (4.3), our studies show that $L_0$ is basically inferior to our test in the sense that it is much less powerful than the proposed test. We now give the corresponding simulated sizes and power values with 1000 replications for model (4.3) for both of the tests in Tables 2–5 below.

Table 1
*Simulated sizes and power values at the 5% level*

| | $T = 250$ | | $T = 500$ | | $T = 750$ | |
|---|---|---|---|---|---|---|
| $\beta$ | $L_0$ | $L_5$ | $L_0$ | $L_5$ | $L_0$ | $L_5$ |
| 0.00 | 0.037 | 0.041 | 0.059 | 0.039 | 0.054 | 0.051 |
| −0.05 | 0.718 | 0.464 | 1.000 | 0.679 | 1.000 | 0.804 |
| −0.10 | 0.999 | 0.811 | 1.000 | 0.966 | 1.000 | 0.986 |
| −0.20 | 1.000 | 0.993 | 1.000 | 1.000 | 1.000 | 1.000 |

The finite-sample results given in Tables 2–5 show that the proposed test and the simulation scheme work well numerically. Table 2 lists the sizes for $L_i$ for $1 \le i \le 5$ and $L_0$. While the sizes are relatively low for $L_5$ in the cases of $T = 250$ and $T = 500$, the size function approaches 5% when $T$ is as large as 750. Most importantly, with such choices of the simulated critical values, Tables 3–5 show that the proposed test is powerful for nonstationarity versus stationarity. For example, when the "distance" between nonstationarity and stationarity is as small as for $\beta = 0.05$, the maximum of the power for $T = 250$ at the 5% level is already over 20%. Comparing the power values of $L_0$ with these values of $L_i$, $1 \le i \le 5$, our observation is that the Dickey–Fuller test is inferior for the case where the alternative is nonlinear. This further supports proposing a test for dealing with such nonparametric nonstationarity.

Table 2
*Simulated sizes at the 5% level*

| $T$ | $L_1$ | $L_2$ | $L_3$ | $L_4$ | $L_5$ | $L_0$ |
|---|---|---|---|---|---|---|
| 250 | 0.003 | 0.010 | 0.034 | 0.047 | 0.039 | 0.038 |
| 500 | 0.007 | 0.017 | 0.026 | 0.041 | 0.037 | 0.061 |
| 750 | 0.005 | 0.014 | 0.038 | 0.050 | 0.049 | 0.056 |

Table 3
*Power values for $T = 250$ at the 5% level*

| $\beta$ | $L_1$ | $L_2$ | $L_3$ | $L_4$ | $L_5$ | $L_0$ |
|---|---|---|---|---|---|---|
| −0.05 | 0.095 | 0.112 | 0.129 | 0.141 | 0.207 | 0.087 |
| −0.10 | 0.206 | 0.268 | 0.350 | 0.438 | 0.647 | 0.127 |
| −0.20 | 0.566 | 0.726 | 0.881 | 0.972 | 0.998 | 0.421 |
| −0.40 | 0.984 | 0.999 | 1.000 | 1.000 | 1.000 | 0.678 |

Table 4
*Power values for $T = 500$ at the 5% level*

| $\beta$ | $L_1$ | $L_2$ | $L_3$ | $L_4$ | $L_5$ | $L_0$ |
|---|---|---|---|---|---|---|
| −0.05 | 0.160 | 0.202 | 0.249 | 0.323 | 0.477 | 0.097 |
| −0.10 | 0.432 | 0.568 | 0.746 | 0.889 | 0.982 | 0.231 |
| −0.20 | 0.923 | 0.993 | 1.000 | 1.000 | 1.000 | 0.519 |
| −0.40 | 1.000 | 1.000 | 1.000 | 1.000 | 1.000 | 0.754 |

TABLE 5
*Power values for $T = 750$ at the 5% level*

| $\beta$ | $L_1$ | $L_2$ | $L_3$ | $L_4$ | $L_5$ | $L_0$ |
|---|---|---|---|---|---|---|
| −0.05 | 0.279 | 0.280 | 0.358 | 0.461 | 0.694 | 0.121 |
| −0.10 | 0.663 | 0.753 | 0.905 | 0.977 | 0.999 | 0.398 |
| −0.20 | 0.992 | 0.999 | 1.000 | 1.000 | 1.000 | 0.689 |
| −0.40 | 1.000 | 1.000 | 1.000 | 1.000 | 1.000 | 0.842 |

As Tables 2–5 show, the corresponding power value of $L_4$ in each case is only the second best among $L_i$ for $1 \le i \le 5$ if we choose an optimal bandwidth such that the simulated size is the closest to 5%. Thus, our finite sample studies also support the fact that there is a kind of trade-off between sizes and power values.

EXAMPLE 4.2. This example examines the three month Treasury Bill rate data given in Figure 1 below sampled monthly over the period from January 1963 to December 1998, providing 432 observations.

Let $\{X_t : t = 1, 2, \ldots, 432\}$ be the set of treasury Bill rate data. As Figure 1 does not suggest that there is any significant trend for the data set, it is not unreasonable to assume that $\{X_t\}$ satisfies a nonlinear autoregressive model of the form

$$X_t = g(X_{t-1}) + e_t \tag{4.5}$$

with the form of $g(\cdot)$ being unknown.

To apply the test $\widehat{L}_T(\widehat{h}_{\text{test}})$ to determine whether $\{X_t\}$ follows a random walk model of the form $X_t = X_{t-1} + u_t$, we need to propose the following procedure for computing the $p$-value of $\widehat{L}_T(\widehat{h}_{\text{test}})$:

- For the real data set, compute $\widehat{h}_{\text{test}}$ and $\widehat{L}_T(\widehat{h}_{\text{test}})$.

- Let $X_1^* = X_1$. Generate a sequence of bootstrap resamples $\{\varepsilon_t^*\}$ from $N(0,1)$ and then $X_t^* = X_{t-1} + \widehat{\sigma}_u \varepsilon_t^*$ for $2 \le t \le 432$.
- Compute the corresponding version $\widehat{L}_T^*(\widehat{h}_{\text{test}})$ of $\widehat{L}_T$ based on $\{X_t^*\}$.
- Repeat the above steps $M$ times to find the bootstrap distribution of $\widehat{L}_T^*(\widehat{h}_{\text{test}})$ and then compute the proportion that $\widehat{L}_T(\widehat{h}_{\text{test}}) < \widehat{L}_T^*(\widehat{h}_{\text{test}})$. This proportion is an approximate $p$-value of $\widehat{L}_T(\widehat{h}_{\text{test}})$.

Our simulation results return the simulated $p$-values of $\widehat{p}_1 = 0.005$ for $L_0$ and $\widehat{p}_2 = 0.011$ for $\widehat{L}_T(\widehat{h}_{\text{test}})$. While both of the simulated $p$-values suggest that there is not enough evidence to accept the unit-root structure at the 5% significance level, there is some evidence of accepting the unit-root structure based on $\widehat{L}_T(\widehat{h}_{\text{test}})$ at the 1% significance level. When we also generated $\{\varepsilon_t^*\}$ from a non-Gaussian distribution, the simulated $p$-values were quite close. By comparison, Jiang (1998) rejects the null hypothesis of nonstationarity on the Fed data based on an application of an augmented Dickey–Fuller unit-root test for $H_0' : \theta = 1$ in a linear model of the form $X_t = \theta X_{t-1} + e_t$.

**5. Conclusion and extensions.** We have proposed a nonparametric specification test for testing whether there is a kind of unit root structure in a nonlinear autoregressive mean function. An asymptotic normal distribution of the proposed test has been established. In addition, we have also proposed a simulation scheme to implement the proposed test in practice. The finite-sample results show that both the proposed test and the simulation scheme are practically applicable and implementable.

It is pointed out that we may also consider a generalized form of model (1.1) with $\sigma_u$ replaced by a stochastic volatility function $\sigma(X_{t-1})$. In this

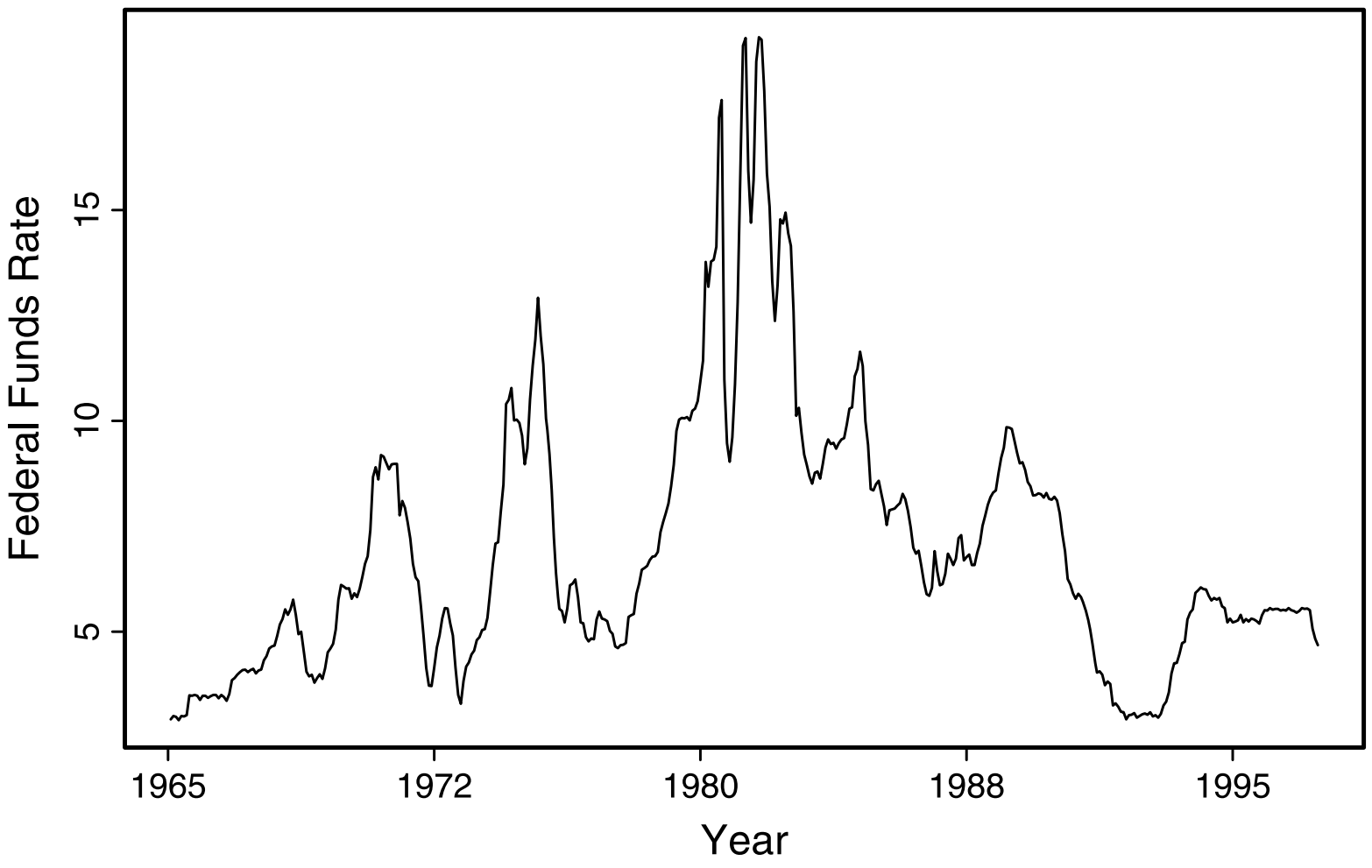

FIG. 1.

case, we should be considering a test for

$$H_{01}: P(g(X_{t-1}) = X_{t-1} \text{ and } \sigma(X_{t-1}) = \sigma_u) = 1. \tag{5.1}$$

In this case, we may use a kernel-based test of the form

$$S_T(h) = \sum_{t=1}^{T} \sum_{s=1, s\neq t}^{T} (U_s K_{h_1}(X_s - X_t) U_t + V_s G_{h_2}(X_s - X_t) V_t), \tag{5.2}$$

where $G_{h_2}(\cdot) = G(\cdot/h_2)$ with $G(\cdot)$ being a probability kernel function, $h = (h_1, h_2)$ is a pair of bandwidth parameters, $U_t = Y_t - \widehat{g}(X_{t-1})$ and $V_t = U_t^2 - \widehat{\sigma}_u^2$ and $\widehat{\sigma}_u$ is an estimator of $\sigma_u$ under $H_0$. Similarly, to Theorems 2.2 and 3.1, we may establish two corresponding theorems for $S_T(h)$. As the details for this case are lengthy and technical, we leave this issue for future study.

Another possible extension will be on the multivariate case where a multivariate autoregressive model is given as follows:

$$X_t = g(X_{t-1}, \ldots, X_{t-p}) + e_t. \tag{5.3}$$

In this case, we are interested in testing a null hypothesis of the form

$$H_{02}: P\left(g(X_{t-1}, \ldots, X_{t-p}) = \sum_{j=1}^{p} \theta_j X_{t-j}\right) = 1, \tag{5.4}$$

in which there is at least one unit root of the corresponding characteristic polynomial. Detailed construction of such a test would involve some estimation procedures for additive models as used in Gao, Lu and Tjøstheim (2006) in the stationary spatial case and as proposed by Gao (2007) in the stationary time series case. Since such an extension is not straightforward, we also leave it as a future topic.

## APPENDIX

This appendix provides the proofs for some necessary technical lemmas that are needed to complete the proofs of Theorems 2.1 and 2.2. Some additional details are given in Appendices B–E of Gao et al. (2008).

Let $a_{st} = K_h(\sum_{i=s}^{t-1} u_i) = K(\frac{\sum_{i=s}^{t-1} u_i}{h})$ and $\eta_t = 2\sum_{s=1}^{t-1} a_{st} u_s$ for $t > s$. Note that we assume without loss of generality that $\sigma_u^2 = E[u_1^2] = 1$ in this appendix.

LEMMA A.1. *Under the conditions of Theorem 2.1, we have under* $H_0$

$$\sigma_T^2 = C_{10} T^{3/2} h (1 + o(1))$$

*for* $T$ *large enough, where* $C_{10} = \frac{16 J_{02}}{3\sqrt{2\pi}}$ *with* $J_{02} = \int K^2(x)\, dx$.

PROOF. Recall $a_{st} = K_h(\sum_{i=s}^{t-1} u_i) = K(\frac{\sum_{i=s}^{t-1} u_i}{h})$ and $\eta_t = 2\sum_{s=1}^{t-1} a_{st} u_s$. It follows under $H_0$ that

$$\begin{aligned}\sigma_T^2 &= \mathrm{E}\left[\sum_{t=2}^{T} \eta_t u_t\right]^2 \\ &= 4\sum_{t=2}^{T}\sum_{s=1}^{t-1} E[a_{st}^2 u_s^2 u_t^2] + 4\sum_{t=2}^{T}\sum_{s_1 \neq s_2 = 1}^{t-1} E[a_{s_1 t} a_{s_2 t} u_{s_1} u_{s_2} u_t^2] \qquad \text{(A.1)} \\ &= 4\sum_{t=2}^{T}\sum_{s=1}^{t-1} E[a_{st}^2 u_s^2] + R_T,\end{aligned}$$

where $R_T = 4\sum_{t=2}^{T}\sum_{s_1\neq s_2=1}^{t-1} E[a_{s_1 t} a_{s_2 t} u_{s_1} u_{s_2}]$.

Let $u_{st} = \sum_{i=s+1}^{t-1} u_i$. Assumption 2.1(i), (ii) already assumes that $\{u_i\}$ is a sequence of independent and identically distributed random variables and has a symmetric probability density function.

Let $f(x)$ and $f_{st}(x)$ be the density functions of $u_i$ and $u_{st}$, respectively, and $g_{st}(x)$ be the density functions of $V_{st} = \frac{u_{st}}{\sqrt{t-s-1}}$. Clearly, $f_{st}(x) = g_{st}(\frac{x}{\sqrt{t-s-1}}) \times \frac{1}{\sqrt{t-s-1}}$, and by utilizing the usual normal approximation of $V_{st} \to_D N(0,1)$ as $t-s\to\infty$ under the conventional central limit theorem conditions, it follows that as $t-s\to\infty$, $g_{st}(x)\to\phi(x)$ and $g_{st}(\frac{x}{\sqrt{t-s-1}}) \to C_0$ uniformly in $x$, where $\phi(x) = \frac{1}{\sqrt{2\pi}}\exp\{-\frac{x^2}{2}\}$, and $C_0 = \phi(0) = \frac{1}{\sqrt{2\pi}}$.

Thus, for $t-s$ large enough, we have

$$\begin{aligned}E[a_{st}^2 u_s^2] &= \int\!\!\int K_h^2(u_{st}+u_s) u_s^2 f(u_s) f_{st}(u_{st})\,du_s\,du_{st} \\ &= h\int\!\!\int K^2(y) x^2 f(x) f_{st}(hy-x)\,dx\,dy \\ &= h(1+o(1))\int\!\!\int K^2(y)x^2 f(x) f_{st}(x)\,dx\,dy \\ &= h(1+o(1))\int\!\!\int K^2(y) x^2 f(x) \\ &\qquad \times g_{st}\left(\frac{x}{\sqrt{t-s-1}}\right)\frac{1}{\sqrt{t-s-1}}\,dx\,dy \qquad \text{(A.2)} \\ &= h(1+o(1))\frac{\int K^2(y)\,dy}{\sqrt{t-s-1}}\int x^2 f(x)\left[g_{st}\left(\frac{x}{\sqrt{t-s-1}}\right)\right]dx \\ &= h(1+o(1))C_0\frac{\int K^2(y)\,dy}{\sqrt{t-s-1}},\end{aligned}$$

where the fact that $\int x^2 f(x)\,dx = E[u_1^2] = 1$ is used.

Choose some positive integer $\Gamma_T \geq 1$ such that $\Gamma_T \to \infty$ and $\frac{\Gamma_T}{\sqrt{T}h} \to 0$ as $T \to \infty$. Observe that

$$\sum_{t=2}^{T}\sum_{s=1}^{t-1} E[a_{st}^2 u_s^2] = \sum_{s=1}^{T-1}\sum_{t=s+1}^{T} E[a_{st}^2 u_s^2] = A_{1T} + A_{2T},$$

where $A_{1T} = \sum_{s=1}^{T-1}\sum_{1\leq(t-s)\leq\Gamma_T} E[a_{st}^2 u_s^2] = O(T\Gamma_T) = o(T^{3/2}h)$ using the fact that $E[a_{st}^2 u_s^2] \leq k_0^2 E[u_s^2] = k_0^2$ due to the boundedness of the kernel $K(\cdot)$ by a constant $k_0 > 0$.

And it follows from (A.2) that

$$\begin{aligned} A_{2T} &= \sum_{s=1}^{T-1}\sum_{\Gamma_T+1\leq(t-s)\leq T-1} E[a_{st}^2 u_s^2] \\ &= \frac{4\int K^2(y)\,dy}{3} C_0 T^{3/2} h(1+o(1)). \end{aligned}$$

It can then be seen that for $T$ large enough

$$\sum_{t=2}^{T}\sum_{s=1}^{t-1} E[a_{st}^2 u_s^2] = \frac{4\int K^2(y)\,dy}{3\sqrt{2\pi}} T^{3/2} h(1+o(1)). \tag{A.3}$$

To deal with $R_T$, we need to introduce the following notation: for $1 \leq i \leq 2$,

$$Z_i = u_{s_i}, \qquad Z_{11} = \sum_{i=s_1+1}^{t-1} u_i, \qquad Z_{22} = \sum_{j=s_2+1}^{s_1-1} u_j, \tag{A.4}$$

ignoring the notational involvement of $s$, $t$ and others.

Let $f_{ii}(x_{ii})$ and $g_{ii}(x_{ii})$ be the probability density functions of $Z_{ii}$ and $\frac{Z_{ii}}{\sigma_{ii}}$, respectively, with $\sigma_{11}^2 = t - s_1 - 1$ and $\sigma_{22}^2 = s_1 - s_2 - 1$.

Clearly, $f_{ii}(x)$ is symmetric due to the symmetry of $f(x)$. Note that $f_{ii}(x) = g_{ii}(\frac{x}{\sigma_{ii}})\frac{1}{\sigma_{ii}}$ and $f'_{ii}(x) = g'_{ii}(\frac{x}{\sigma_{ii}})\frac{1}{\sigma_{ii}^2}$.

By utilizing the normal approximation of $\frac{Z_{ii}}{\sigma_{ii}} \to_D N(0,1)$ as $\sigma_{ii} \to \infty$ under the usual central limit theorem conditions, it follows that $g_{ii}(x) \to \phi(x)$ and $g_{ii}(\frac{x}{\sigma_{ii}}) \to C_0$, with $C_0 = \frac{1}{\sqrt{2\pi}}$, and $\frac{1}{x}g'_{ii}(x) \to \frac{1}{x}\phi'(x) = -\phi(x)$ leading to $\frac{\sigma_{ii}}{x}g'_{ii}(\frac{x}{\sigma_{ii}}) \to -\phi(0) = -C_0$, as $\sigma_{ii} \to \infty$.

Similarly to (A.2), we can derive that as $\sigma_{11}^2 = t - s_1 - 1 \to \infty$ and $\sigma_{22}^2 = s_1 - s_2 - 1 \to \infty$,

$$\begin{aligned} &E[a_{s_1t}a_{s_2t}u_{s_1}u_{s_2}] \\ &\quad = E\left[K_h\left(\sum_{i=s_1}^{t-1} u_i\right)K_h\left(\sum_{j=s_2}^{t-1} u_j\right)u_{s_1}u_{s_2}\right] \\ &\quad = E[Z_1Z_2K_h(Z_1+Z_{11})K_h(Z_1+Z_2+Z_{11}+Z_{22})] \end{aligned}$$

$$
\begin{aligned}
&= E\left[\prod_{i=1}^{2} Z_i K_h\left(\sum_{j=1}^{i}(Z_j + Z_{jj})\right)\right] \\
&= \int\cdots\int x_1 x_2 K_h(x_1 + x_2 + x_{11} + x_{22}) K_h(x_1 + x_{11}) \\
&\qquad \times f(x_1) f(x_2) f_{11}(x_{11}) f_{22}(x_{22})\, dx_1\, dx_2\, dx_{11}\, dx_{22} \\
&\qquad\qquad \text{using } y_{ii} = \frac{x_i + x_{ii}}{h} \\
&= h^2 \prod_{j=1}^{2}\left[\left(\int\int K\left(\sum_{i=1}^{j} y_{ii}\right) x_j f(x_j) f_{jj}(x_j - h y_{jj})\, dx_j\, dy_{jj}\right)\right] \\
&\qquad \left(\text{using Taylor expansions and } \int x_j f(x_j) f_{jj}(x_j)\, dx_j = 0 \right. \\
&\qquad\qquad \left.\text{due to symmetry of } f \text{ and } g_{jj}\right) \\
&= h^4(1 + o(1)) \prod_{j=1}^{2}\left[\int\int y_{jj} K\left(\sum_{i=1}^{j} y_{ii}\right) x_j f(x_j) f'_{jj}(x_j)\, dx_j\, dy_{jj}\right] \\
&= h^4(1 + o(1)) \prod_{j=1}^{2}\left[\int y_{jj} K\left(\sum_{i=1}^{j} y_{ii}\right) dy_{jj} \cdot \int x_j f(x_j) f'_{jj}(x_j)\, dx_j\right] \\
&= \frac{C_{11}(K) h^4 (1 + o(1))}{\prod_{j=1}^{2} \sigma_{jj}^2} \prod_{j=1}^{2}\left[\int x_j f(x_j) g'_{jj}\left(\frac{x_j}{\sigma_{jj}}\right) dx_j\right] \\
&= \frac{C_{11}(K) h^4 (1 + o(1))}{\prod_{j=1}^{2} \sigma_{jj}^2} \prod_{j=1}^{2}\left[\int \frac{x_j^2}{\sigma_{jj}} f(x_j) \frac{\sigma_{jj}}{x_j} g'_{jj}\left(\frac{x_j}{\sigma_{jj}}\right) dx_j\right] \\
&= \frac{C_{11}(K) h^4 (1 + o(1))}{2\pi} \prod_{j=1}^{2} \frac{1}{(\sqrt{1 + \sigma_{jj}^2})^3} \\
&= O(h^4) \frac{1}{(\sqrt{t - s_1})^3} \frac{1}{(\sqrt{s_1 - s_2})^3},
\end{aligned} \tag{A.5}
$$

where the conventional notation $\prod$ defines $\prod_{i=1}^{k} p_i = p_1 p_2 \cdots p_k$, and

$$
\begin{aligned}
C_{11}(K) &= \prod_{j=1}^{2}\left(\int y_{jj} K\left(\sum_{i=1}^{j} y_{ii}\right) dy_{jj}\right) \\
&= \int\int y_{11} y_{22} K(y_{11}) K(y_{11} + y_{22})\, dy_{11}\, dy_{22} = -\int y^2 K(y)\, dy.
\end{aligned}
$$

Choose $\Gamma_T$ satisfying $\Gamma_T \to \infty$ and $\frac{\Gamma_T^2}{\sqrt{T}h} \to 0$ as $T \to \infty$. Note that

$$\sum_{t=2}^{T} \sum_{s_1 \neq s_2 = 1}^{t-1} E[a_{s_1 t} a_{s_2 t} u_{s_1} u_{s_2}] = A_{3T} + A_{4T} + A_{5T} + A_{6T},$$

where $A_{3T} = \sum_{s_2=1}^{T-2} \sum_{s_1=s_2+1}^{s_2+\Gamma_T} \sum_{t=s_1+1}^{s_1+\Gamma_T} E[a_{s_1 t} a_{s_2 t} u_{s_1} u_{s_2}] = O(T\Gamma_T^2) = o(T^{3/2}h)$ owing to $E[a_{s_1 t} a_{s_2 t} u_{s_1} u_{s_2}] \leq k_0^2 E|u_{s_1} u_{s_2}| \leq k_0^2$ by the assumption that $K(\cdot)$ is bounded by $k_0$:

$$\begin{aligned}
A_{4T} &= \sum_{s_2=1}^{T-2} \sum_{s_1=s_2+1}^{s_2+\Gamma_T} \sum_{t=s_1+\Gamma_T+1}^{T} E[a_{s_1 t} a_{s_2 t} u_{s_1} u_{s_2}] \\
&\leq \sum_{s_2=1}^{T-2} \sum_{s_1=s_2+1}^{s_2+\Gamma_T} \sum_{t=s_1+\Gamma_T+1}^{T} (E[a_{s_1 t}^2 u_{s_1}^2])^{1/2} (E[a_{s_2 t}^2 u_{s_2}^2])^{1/2} \\
&= O(1) \sum_{s_2=1}^{T-2} \sum_{s_1=s_2+1}^{s_2+\Gamma_T} \sum_{t=s_1+\Gamma_T+1}^{T} [h(t-s_1-1)^{1/2}]^{1/2} \\
&= O(\Gamma_T T^{1+1/4} h^{1/2}) = o(T^{3/2}h).
\end{aligned}$$

Similarly to $A_{4T}$, we have

$$A_{5T} = \sum_{s_2=1}^{T-2} \sum_{s_1=s_2+\Gamma_T+1}^{T} \sum_{t=s_1+1}^{t=s_1+\Gamma_T} E[a_{s_1 t} a_{s_2 t} u_{s_1} u_{s_2}] = o(T^{3/2}h).$$

Finally, owing to (A.5),

$$\begin{aligned}
A_{5T} &= \sum_{s_2=1}^{T-2} \sum_{s_1=s_2+\Gamma_T+1}^{T} \sum_{t=s_1+\Gamma_T+1}^{t=T} E[a_{s_1 t} a_{s_2 t} u_{s_1} u_{s_2}] \\
&= O(h^4) \sum_{s_2=1}^{T-2} \sum_{s_1=s_2+\Gamma_T+1}^{T} \sum_{t=s_1+\Gamma_T+1}^{T-1} \frac{1}{(\sqrt{t-s_1})^3} \frac{1}{(\sqrt{s_1-s_2})^3} \\
&= O(T^2 h^4) = o(T^{3/2}h).
\end{aligned}$$

Thus, for $T$ large enough

$$\begin{aligned}
\sum_{t=2}^{T} \sum_{s_1 \neq s_2 = 1}^{t-1} E[a_{s_1 t} a_{s_2 t} u_{s_1} u_{s_2}] &= 2 \sum_{s_2=1}^{T-2} \sum_{s_1=s_2+1}^{T-1} \sum_{t=s_1+1}^{T} E[a_{s_1 t} a_{s_2 t} u_{s_1} u_{s_2}] \\
&= o(T^{3/2}h)
\end{aligned} \tag{A.6}$$

using Assumption 2.1.

Therefore, (A.3) and (A.6) show that as $T \to \infty$

$$\sigma_T^2 = \frac{16 \int K^2(y)\,dy}{3\sqrt{2\pi}} T^{3/2} h(1 + o(1)). \tag{A.7}$$

The proof of Lemma A.1 is therefore finished. $\square$

LEMMA A.2. *Assume that the probability space $(\Omega_n, \mathcal{F}_n, P_n)$ supports square integrable random variables $S_{n,1}, S_{n,2}, \ldots, S_{n,k_n}$, and that the $S_{n,t}$ are adapted to $\sigma$-algebras $\mathcal{F}_{n,t}$, $1 \le t \le k_n$, where*

$$\mathcal{F}_{n,1} \subset \mathcal{F}_{n,2} \subset \cdots \subset \mathcal{F}_{n,k_n} \subset \mathcal{F}_n.$$

*Let $X_{n,t} = S_{n,t} - S_{n,t-1}$, $S_{n,0} = 0$ and $U_{n,t}^2 = \sum_{s=1}^t X_{n,s}^2$. If $\mathcal{G}_n$ is a sub-$\sigma$-algebra of $\mathcal{F}_n$, let $\mathcal{G}_{n,t} = \mathcal{F}_{n,t} \vee \mathcal{G}_n$ (the $\sigma$-algebra generated by $\mathcal{F}_{n,t} \cup \mathcal{G}_n$) and let $\mathcal{G}_{n,0} = \{\Omega_n, \phi\}$ denote the trivial $\sigma$-algebra. Moreover, suppose that*

$$\sum_{t=1}^n E(X_{n,t}^2 I[|X_{n,t}| > \delta] | \mathcal{G}_{n,t-1}) \to_P 0 \tag{A.8}$$

*for some $\delta > 0$, and there exists a $\mathcal{G}_n$-measurable random variable $u_n^2$, such that*

$$U_{n,k_n}^2 - u_n^2 \to_P 0, \tag{A.9}$$

$$\sum_{t=1}^n E(X_{n,t} | \mathcal{G}_{n,t-1}) \to_P 0, \tag{A.10}$$

$$\sum_{t=1}^n |E(X_{n,t} | \mathcal{G}_{n,t-1})|^2 \to_P 0. \tag{A.11}$$

*If*

$$\lim_{\delta \to 0} \lim_{n \to \infty} \inf P\{U_{n,k_n} > \delta\} = 1, \tag{A.12}$$

*then as $n \to \infty$*

$$\frac{S_{n,k_n}}{U_{n,k_n}} \to_D N(0,1).$$

PROOF. The proof follows from Corollary 3.1 and Theorem 3.4 of Hall and Heyde (1980). $\square$

LEMMA A.3. *Under the conditions of Theorem 2.2, we have as $T \to \infty$*

$$\frac{M_{T1}}{\widetilde{\sigma}_T} \to_D N(0,1). \tag{A.13}$$

PROOF. We apply Lemma A.2 to prove Lemma A.3. Let $Y_{Tt} = \frac{\eta_t u_t}{\sigma_T}$, $\Omega_{T,s} = \sigma\{Y_{Tt} : 1 \le t \le s\}$ be a $\sigma$-field generated by $\{Y_{Tt} : 1 \le t \le s\}$, $\mathcal{G}_T = \Omega_{T,P(T)}$ and $\mathcal{G}_{T,s}$ be defined by

$$\mathcal{G}_{T,s} = \begin{cases} \Omega_{T,P(T)}, & 1 \le s \le P(T), \\ \Omega_{T,s}, & P(T)+1 \le s \le T, \end{cases} \tag{A.14}$$

where $P(T) \ge 1$ is chosen such that $P(T) \to \infty$ and $\frac{P(T)}{T} \to 0$ as $T \to \infty$. Let $\widetilde{U}^2_{P(T)} = \frac{\widetilde{\sigma}^2_{P(T)}}{\sigma^2_{P(T)}}$, where $\widetilde{\sigma}^2_S = 2\sum_{t=1}^S \sum_{s=1,s\ne t}^S u_s^2 a_{st}^2 u_t^2$ and $\sigma^2_S = \mathrm{var}[\sum_{t=2}^S \eta_t u_t]$ for all $1 \le S \le T$ as defined before.

In view of Lemma A.2 above, in order to prove that as $T \to \infty$

$$\frac{M_{T1}}{\widetilde{\sigma}_T} = \frac{1}{\widetilde{\sigma}_T}\sum_{t=2}^T \eta_t u_t \to_D N(0,1), \tag{A.15}$$

it suffices to show that for all $\delta > 0$,

$$\sum_{t=2}^T E[Y_{Tt}^2 I_{\{[|Y_{Tt}|>\delta]\}} | \Omega_{T,t-1}] \to_P 0, \tag{A.16}$$

$$\frac{\widetilde{\sigma}^2_T}{\sigma^2_T} - \widetilde{U}^2_{P(T)} \to_P 0, \tag{A.17}$$

$$\begin{aligned} \sum_{t=2}^T E[Y_{Tt}|\mathcal{G}_{T,t-1}] &= \sum_{t=2}^{P(T)} Y_{Tt} + \sum_{t=P(T)+1}^T E[Y_{Tt}|\Omega_{T,t-1}] \\ &= \sum_{t=2}^{P(T)} Y_{Tt} \to_P 0, \end{aligned} \tag{A.18}$$

$$\begin{aligned} \sum_{t=2}^T |E[Y_{Tt}|\mathcal{G}_{T,t-1}]|^2 &= \sum_{t=2}^{P(T)} Y_{Tt}^2 + \sum_{t=P(T)+1}^T |E[Y_{Tt}|\Omega_{T,t-1}]|^2 \\ &= \sum_{t=2}^{P(T)} Y_{Tt}^2 \to_P 0, \end{aligned} \tag{A.19}$$

$$\lim_{\delta \to 0} \lim \inf_{T \to \infty} P\left(\frac{\widetilde{\sigma}_T}{\sigma_T} > \delta\right) = 1. \tag{A.20}$$

The proof of (A.18) is similar to that of (A.19), which follows from

$$\sum_{t=2}^{P(T)} E[Y_{Tt}^2] = O\left(\left(\frac{P(T)}{T}\right)^{3/2}\right) \to 0 \tag{A.21}$$

as $T \to \infty$, in which Lemma A.1 has been used.

In order to prove (A.16), it suffices to show that

$$\frac{1}{\sigma_T^4}\sum_{t=2}^{T} E[\eta_t^4] \to 0. \tag{A.22}$$

The proof of (A.22) is given in Lemma A.4 below. The proof of (A.17) is given in Lemma A.5 below.

The proof of (A.20) follows from

$$\frac{\widetilde{\sigma}_T^2}{\sigma_T^2} \to_D \xi^2 > 0 \tag{A.23}$$

for some random variable $\xi^2$. The proof is given in Lemma A.6 below. □

Lemma A.4. *Under the conditions of Theorem 2.2, we have*

$$\lim_{T\to\infty}\frac{1}{\sigma_T^4}\sum_{t=2}^{T} E[\eta_t^4] = 0. \tag{A.24}$$

Proof. Observe that

$$E[\eta_t^4] = 16\sum_{s_1=1}^{t-1}\sum_{s_2=1}^{t-1}\sum_{s_3=1}^{t-1}\sum_{s_4=1}^{t-1} E[a_{s_1t}a_{s_2t}a_{s_3t}a_{s_4t}u_{s_1}u_{s_2}u_{s_3}u_{s_4}]. \tag{A.25}$$

We mainly consider the cases of $s_i \neq s_j$ for all $i \neq j$ in the following proof. Since the other terms involve at most triple summations, we may deal with such terms similarly. Without loss of generality, we only look at the case of $1 \le s_4 < s_3 < s_2 < s_1 \le t-1$ in the following evaluation. Let

$$\sum_{i=s_1}^{t-1} u_i = u_{s_1} + \sum_{i=s_1+1}^{t-1} u_i,$$

$$\sum_{i=s_2}^{t-1} u_i = u_{s_1} + u_{s_2} + \sum_{i=s_2+1}^{s_1-1} u_i + \sum_{j=s_1+1}^{t-1} u_j,$$

$$\sum_{i=s_3}^{t-1} u_i = u_{s_1} + u_{s_2} + u_{s_3} + \sum_{k=s_3+1}^{s_2-1} u_k + \sum_{i=s_2+1}^{s_1-1} u_i + \sum_{j=s_1+1}^{t-1} u_j,$$

$$\sum_{i=s_4}^{t-1} u_i = u_{s_1} + u_{s_2} + u_{s_3} + u_{s_4} + \sum_{l=s_4+1}^{s_3-1} u_l + \sum_{k=s_3+1}^{s_2-1} u_k + \sum_{i=s_2+1}^{s_1-1} u_i + \sum_{j=s_1+1}^{t-1} u_j.$$

Similarly to (A.4), let again $Z_i = u_{s_i}$ for $1 \le i \le 4$,

$$Z_{11} = \sum_{i=s_1+1}^{t-1} u_i, \qquad Z_{22} = \sum_{j=s_2+1}^{s_1-1} u_j,$$

$$Z_{33} = \sum_{k=s_3+1}^{s_2-1} u_k, \qquad Z_{44} = \sum_{l=s_4+1}^{s_3-1} u_l.$$

By the same arguments as in the proof of (A.5), we have

$$\begin{aligned}
E\Bigg[\prod_{i=1}^{4} a_{s_i t} u_{s_i}\Bigg] &= E\Bigg[\prod_{j=1}^{4} Z_j K_h\Bigg(\sum_{i=1}^{j}[Z_i + Z_{ii}]\Bigg)\Bigg] \\
&= \prod_{j=1}^{4}\Bigg(\int K_h\Bigg(\sum_{i=1}^{j}[x_i + x_{ii}]\Bigg) x_j f(x_j) f_{jj}(x_{jj})\,dx_j\,dx_{jj}\Bigg) \\
&\qquad\qquad \text{using } y_{ii} = \frac{x_i + x_{ii}}{h} \\
&\quad + h^4 \prod_{j=1}^{4}\Bigg[\int K\Bigg(\sum_{i=1}^{j} y_{ii}\Bigg) x_j f(x_j) f_{jj}(x_j - h y_{jj})\,dx_j\,dy_{jj}\Bigg] \\
&\qquad\qquad \text{using Taylor expansions and } \int x_j f(x_j) f_{jj}(x_j)\,dx_j = 0 \\
&= h^8(1+o(1)) \prod_{j=1}^{4}\Bigg[\int y_{jj} K(u_{jj}) x_j f(x_j) f'_{jj}(x_j)\,dx_j\,dy_{jj}\Bigg] \\
&= h^8(1+o(1)) \prod_{j=1}^{4}\Bigg[\int y_{jj} K(u_{jj})\,dy_{jj} \cdot \int x_j f(x_j) f'_{jj}(x_j)\,dx_j\Bigg] \\
&\qquad\qquad \text{using } f'_{ii}(x) = g'_{ii}\bigg(\frac{x}{\sigma_{ii}}\bigg)\frac{1}{\sigma_{ii}^2} \\
&= \frac{C_{22}(K) h^8 (1+o(1))}{\prod_{j=1}^{4} \sigma_{jj}^2} \prod_{j=1}^{4}\Bigg[\int \frac{x_j^2}{\sigma_{jj}} f(x_j) \frac{\sigma_{jj}}{x_j} g'_{jj}\bigg(\frac{x_j}{\sigma_{jj}}\bigg)\,dx_j\Bigg] \\
&= \frac{C_{22}(K) h^8 (1+o(1))}{4\pi^2} \prod_{j=1}^{4} \frac{1}{(1+\sigma_{jj}^2)\sigma_{jj}},
\end{aligned} \tag{A.26}$$

where $u_{jj} = \sum_{i=1}^{j} y_{ii}$ is used to shorten some expressions, and

$$C_{22}(K) = \prod_{j=1}^{4}\Bigg(\int y_{jj} K\Bigg(\sum_{i=1}^{j} y_{ii}\Bigg)\,dy_{jj}\Bigg) < \infty.$$

Hence, similarly to (A.3), we have

$$\sum_{t=2}^{T} \sum_{1 \le s_4 < s_3 < s_2 < s_1 \le t-1} E[a_{s_1 t} a_{s_2 t} a_{s_3 t} a_{s_4 t} u_{s_1} u_{s_2} u_{s_3} u_{s_4}] = o(T^3 h^2) \tag{A.27}$$

using Assumption 2.1.

Analogously, we can deal with the other terms of (A.25) as follows:

$$\sum_{t=2}^{T}\sum_{1\le s_2\ne s_1\le t-1} E[a_{s_1t}^2 a_{s_2t}^2 u_{s_1}^2 u_{s_2}^2] = o(T^3h^2), \tag{A.28}$$

$$\sum_{t=2}^{T}\sum_{1\le s_3\ne s_2\ne s_1\le t-1} E[a_{s_1t}^2 a_{s_2t} a_{s_3t} u_{s_1}^2 u_{s_2} u_{s_3}] = o(T^3h^2), \tag{A.29}$$

$$\sum_{t=2}^{T}\sum_{1\le s_2\ne s_1\le t-1} E[a_{s_1t}^3 a_{s_2t} u_{s_1}^3 u_{s_2}] = o(T^3h^2). \tag{A.30}$$

Thus, we have finished the proof of (A.24) using (A.25)–(A.30). □

Lemma A.5. *Let the conditions of Theorem 2.2 hold. Then as* $T\to\infty$

$$\frac{\widetilde{\sigma}_T^2}{\sigma_T^2} - \widetilde{U}_{P(T)}^2 \to_P 0. \tag{A.31}$$

Proof. For $1\le S\le T$, recall $\widetilde{U}_S^2 = \frac{\widetilde{\sigma}_S^2}{\sigma_S^2}$, where $\widetilde{\sigma}_S^2 = 2\sum_{t=1}^{S}\sum_{s=1,s\ne t}^{S} u_s^2 \times a_{st}^2 u_t^2$.

To use simplified notation in the proof of this lemma, we introduce the following lower-case notation: $m=T$, $n=P(T)$, $\sigma_m^2=\sigma_T^2$, $\sigma_n^2=\sigma_{P(T)}^2$, and for $1\le i\le n$, $1\le j\le i-1$,

$$e_{ij} = (u_i^2 - E[u_1^2])K_h^2\left(\sum_{l=j}^{i-1} u_l\right)u_j^2 \quad\text{and}\quad X_{mi} = \frac{1}{\sigma_m^2}\sum_{j=1}^{i-1} e_{ij}, \tag{A.32}$$

$$v_i^2 = \sum_{j=1}^{i-1} K_h^2\left(\sum_{l=j}^{i-1} u_l\right)u_j^2 = \sum_{j=1}^{i-1} K_h^2\left(\sum_{l=j+1}^{i-1} u_l + u_j\right)u_j^2. \tag{A.33}$$

Note that $X_{mi} = \frac{1}{\sigma_m^2}(u_i^2 - E[u_1^2])v_i^2$ with $E[X_{mi}]=0$.

Observe that

$$\begin{aligned}\frac{\widetilde{\sigma}_m^2}{\sigma_m^2} - \frac{\widetilde{\sigma}_n^2}{\sigma_n^2} &= \sum_{i=1}^{m} X_{mi} - \sum_{j=1}^{n} X_{nj} + E[u_i^2]\left(\frac{1}{\sigma_m^2}\sum_{i=1}^{m} v_i^2 - \frac{1}{\sigma_n^2}\sum_{j=1}^{n} v_j^2\right)\\ &\equiv I_{mn} + E[u_1^2]J_{mn},\end{aligned} \tag{A.34}$$

where $I_{mn} = \sum_{i=1}^{m} X_{mi} - \sum_{j=1}^{n} X_{nj}$ and $J_{mn} = \frac{1}{\sigma_m^2}\sum_{i=1}^{m} v_i^2 - \frac{1}{\sigma_n^2}\sum_{j=1}^{n} v_j^2$.

In view of (A.33), in order to prove (A.31), it suffices to show that as $m,n\to\infty$

$$I_{mn} \to_P 0 \quad\text{and}\quad J_{mn}\to_P 0. \tag{A.35}$$

We now prove the first part of (A.35). In view of the fact that the independence of $\{u_i\}$ implies for $n+1 \le i \le m$ and $1 \le j \le n$,

$$E[X_{mi}(X_{mj}-X_{nj})] = \frac{\sigma_n^2-\sigma_m^2}{\sigma_m^4\sigma_n^2}\sum_{k=1}^{i-1}\sum_{l=1}^{j-1}E[(u_i^2-E[u_1^2])] \times E\left[(u_j^2-E[u_1^2])K_h^2\left(\sum_{p=k}^{i-1}u_p\right)u_k^2K_h^2\left(\sum_{q=l}^{j-1}u_q\right)u_l^2\right]=0,$$

we have

$$\begin{aligned} E[I_{mn}^2] &= E\left[\sum_{i=1}^m X_{mi}-\sum_{j=1}^n X_{nj}\right]^2 = E\left[\sum_{i=n+1}^m X_{mi}+\sum_{j=1}^n (X_{mj}-X_{nj})\right]^2 \\ &= E\left[\sum_{i=n+1}^m X_{mi}\right]^2 + E\left[\sum_{j=1}^n (X_{mj}-X_{nj})\right]^2 \\ &= \frac{1}{\sigma_m^4}\sum_{i=n+1}^m E(u_i^2-E[u_i^2])^2E[v_i^4]+\frac{(\sigma_m^2-\sigma_n^2)^2}{\sigma_m^4\sigma_n^4} \\ &\quad \times \sum_{j=1}^n E(u_j^2-E[u_1^2])^2E[v_j^4]. \end{aligned} \tag{A.36}$$

We start by looking at $\sum_{i=n+1}^m E[v_i^4]$ and $\sum_{j=1}^n E[v_j^4]$ in order to complete the proof of the first part of (A.35). Before we compute the two terms, we have a look at how to prove the second part of (A.35). Note that

$$\begin{aligned} E[J_{mn}^2] &= E\left[\frac{1}{\sigma_m^2}\sum_{i=1}^m v_i^2-\frac{1}{\sigma_n^2}\sum_{j=1}^n v_j^2\right]^2 \\ &= E\left[\frac{1}{\sigma_m^2}\sum_{i=n+1}^m v_i^2+\frac{\sigma_n^2-\sigma_m^2}{\sigma_m^2\sigma_n^2}\sum_{j=1}^n v_j^2\right]^2 \\ &= \frac{1}{\sigma_m^4}E\left[\sum_{i=n+1}^m v_i^2\right]^2+\frac{(\sigma_n^2-\sigma_m^2)^2}{\sigma_m^4\sigma_n^4}E\left[\sum_{j=1}^n v_j^2\right]^2 \\ &\quad + 2\frac{\sigma_n^2-\sigma_m^2}{\sigma_m^4\sigma_n^2}\sum_{i=n+1}^m\sum_{j=1}^n E[v_i^2v_j^2]. \end{aligned} \tag{A.37}$$

We first deal with (A.37) term by term. Recalling $a_{ji} = K_h(\sum_{l=j}^{i-1} u_l)$, we have

$$E\left(\sum_{i=n+1}^{m} v_i^2\right)^2 = E\left[\sum_{i=n+1}^{m}\sum_{j=n+1}^{m} v_i^2 v_j^2\right]$$
$$= \sum_{i=n+1}^{m} E[v_i^4] + \sum_{i=n+1}^{m}\sum_{j=n+1, j\neq i}^{m} E[v_i^2 v_j^2]. \tag{A.38}$$

Observe that

$$E[v_j^2 v_i^2] = \sum_{c=1}^{i-1}\sum_{d=1}^{j-1} E[a_{ci}^2 u_c^2 a_{dj}^2 u_d^2]$$
$$= \sum_{c=1}^{j-1}\sum_{d=1}^{j-1} E[a_{ci}^2 u_c^2 a_{dj}^2 u_d^2] + \sum_{c=j}^{i-1}\sum_{d=1}^{j-1} E[a_{ci}^2 u_c^2 a_{dj}^2 u_d^2] \tag{A.39}$$
$$\equiv I_{ij} + J_{ij},$$

where

$$I_{ij} = \sum_{c=1}^{j-1}\sum_{d=1}^{j-1} E[a_{ci}^2 u_c^2 a_{dj}^2 u_d^2]$$
$$= \sum_{c=1}^{j-1} E[a_{ci}^2 a_{cj}^2 u_c^4] + 2\sum_{c=2}^{j-1}\sum_{d=1}^{c-1} E[a_{ci}^2 u_c^2 a_{dj}^2 u_d^2] \tag{A.40}$$
$$\equiv I_{ij}(1) + I_{ij}(2),$$

$$J_{ij} = \sum_{c=j}^{i-1}\sum_{d=1}^{j-1} E[a_{ci}^2 u_c^2 a_{dj}^2 u_d^2] = \sum_{c=j}^{i-1}\sum_{d=1}^{j-1} E[a_{ci}^2 u_c^2] E[a_{dj}^2 u_d^2] \tag{A.41}$$

using the fact that $\{u_k : j \le k \le i-1\}$ and $\{u_l : 1 \le l \le j-1\}$ are all mutually independent.

Thus, we need only to evaluate $\sum_{i=2}^{n}\sum_{j=1}^{i-1} I_{ij}$. To do so, we introduce another set of simplified symbols: $Z_{11} = \sum_{k=d+1}^{c-1} u_k$, $Z_{22} = \sum_{k=c+1}^{j-1} u_k$, $Z_{33} = \sum_{l=j}^{i-1} u_l$, $Z_1 = u_d$ and $Z_2 = u_c$. In this case, we have the following decompositions: for $1 \le d \le c-1$, $1 \le d \le j-1$ and $1 \le j \le i-1$,

$$\sum_{l=c}^{i-1} u_l = u_c + \sum_{l=c+1}^{j-1} u_l + \sum_{l=j}^{i-1} u_l = Z_2 + Z_{22} + Z_{33},$$

$$\sum_{k=d}^{j-1} u_k = u_d + \sum_{k=d+1}^{c-1} u_k + u_c + \sum_{k=c+1}^{j-1} u_k = Z_1 + Z_2 + Z_{11} + Z_{22}.$$

By the same arguments as used in the proof of Lemma A.1, we have

$$
\begin{aligned}
E[a_{ci}^2 u_c^2 a_{cj}^2 u_c^2] &= E[K_h^2(Z_2+Z_{22})K_h^2(Z_2+Z_{22}+Z_{33})Z_2^4] \\
&= \int\cdots\int K_h^2(x_2+x_{22})K_h^2(x_2+x_{22}+x_{33}) \\
&\qquad \times x_2^4 f(x_2) f_{22}(x_{22}) f_{33}(x_{33})\,dx_2\,dx_{22}\,dx_{33} \\
&\qquad\qquad \text{using } y_2 = x_2, y_{22} = \frac{x_2+x_{22}}{h}, y_{33} = \frac{x_{33}}{h} \\
&= h^2 \int\cdots\int K^2(y_{22})K^2(y_{22}+y_{33})y_2^4 \\
&\qquad \times f(y_2) f_{22}(y_2 - y_{22}h) f_{33}(hy_{33})\,dy_2\,dy_{22}\,dy_{33} \\
&= h^2(1+o(1))\left(\int K^2(u)\,du\right)^2 \\
&\quad \times \left(\int x_2^4 f(x_2) f_{22}(x_2)\,dx_2\right) f_{33}(0),
\end{aligned}
$$

where $f_{ii}(\cdot)$ denotes the marginal density of $Z_{ii}$ and $f(\cdot)$ denotes the density of $Z_i$.

Similarly, we have

$$
\begin{aligned}
E[a_{ci}^2 u_c^2 a_{dj}^2 u_d^2] &= E\left[K_h^2\left(\sum_{i=1}^2 (Z_i+Z_{ii})\right) K_h^2(Z_2+Z_{22}+Z_{33})Z_1^2 Z_2^2\right] \\
&= \int\cdots\int K_h^2\left(\sum_{i=1}^2 (x_i+x_{ii})\right) K_h^2(x_2+x_{22}+x_{33}) \\
&\qquad \times \left(\prod_{i=1}^2 x_i^2 f(x_i) f_{ii}(x_{ii})\,dx_i\,dx_{ii}\right) f_{33}(x_{33})\,dx_{33} \\
&= h^3(1+o(1)) \int\cdots\int K^2(y_{11}+y_{22})K^2(y_{22}+y_{33})y_1^2 y_2^2 \\
&\qquad \times f(y_1) f(y_2) f_{11}(y_1 - y_{11}h) \\
&\qquad \times f_{22}(y_2 - y_{22}h) f_{33}(0) \\
&\qquad \times dy_1\,dy_2\,dy_{11}\,dy_{22}\,dy_{33} \\
&= h^3(1+o(1))\left(\int K^2(u)\,du\right)^2 f_{33}(0) \\
&\quad \times \left(\int x_1^2 f(x_1) f_{11}(x_1)\,dx_1\right)\left(\int x_2^2 f(x_2) f_{22}(x_2)\,dx_2\right).
\end{aligned}
$$

Using the same arguments as used in the calculations of (A.2), (A.3) and (A.7), we have

$$\begin{aligned}32\sum_{i=n+1}^{m}\sum_{j=1}^{i-1}I_{ij}(1) &= C_{10}^2(m-n)^3h^2(1+o(1))\\ &= \sigma_{m-n}^4(1+o(1)),\end{aligned}\tag{A.42}$$

$$\begin{aligned}32\sum_{i=n+1}^{m}\sum_{j=1}^{i-1}I_{ij}(2) &= C(m-n)^2h^3(1+o(1))\\ &= o(\sigma_{m-n}^4),\end{aligned}\tag{A.43}$$

where $C_{10}$ is as defined in Theorem 2.1 and $C > 0$ is a positive constant.

Similarly, by (A.41) we have

$$\begin{aligned}32\sum_{i=n+1}^{m}\sum_{j=1}^{i-1}J_{ij}(2) &= 32\sum_{i=n+1}^{m}\sum_{j=1}^{i-1}\sum_{c=j}^{i-1}\sum_{d=1}^{j-1}E[a_{ci}^2u_c^2]E[a_{dj}^2u_d^2]\\ &= o(\sigma_{m-n}^4).\end{aligned}\tag{A.44}$$

Hence, (A.39)–(A.44) imply for $m$ and $n$ large enough,

$$\begin{aligned}E\left[\sum_{i=n+1}^{m}\sum_{j=n+1,j\neq i}^{m}v_i^2v_j^2\right] &= \sum_{i=n+1}^{m}\sum_{j=n+1,j\neq i}^{m}E[v_i^2v_j^2]\\ &= \sigma_{m-n}^4(1+o(1)),\end{aligned}\tag{A.45}$$

where $\sigma_m^2$ is as defined above (A.32).

Analogously to (A.5), we can show that for $m$ and $n$ large enough,

$$\begin{aligned}\sum_{i=n+1}^{m}E[v_i^4] &= \sum_{i=n+1}^{m}\sum_{s=1}^{i-1}\sum_{t=1}^{i-1}E[a_{si}^2a_{ti}^2u_s^2u_t^2]\\ &= O(h^2)\sum_{i=n+1}^{m}\sum_{s=2}^{i-1}\sum_{t=1}^{s-1}\frac{1}{\sqrt{i-s}}\frac{1}{\sqrt{s-t}}\\ &= O(h^2(m^2-n^2)) = o(\sigma_{m-n}^4).\end{aligned}\tag{A.46}$$

Similarly to (A.45), we can show that for $n$ large enough,

$$\begin{aligned}\sum_{i=1}^{n}E[v_i^4] &= \sum_{i=1}^{n}\sum_{s=1}^{i-1}\sum_{t=1}^{i-1}E[a_{si}^2a_{ti}^2u_s^2u_t^2]\\ &= O(h^2)\sum_{i=1}^{n}\sum_{s=2}^{i-1}\sum_{t=1}^{s-1}\frac{1}{\sqrt{i-s}}\frac{1}{\sqrt{s-t}}\\ &= O(h^2n^2) = o(\sigma_n^4).\end{aligned}\tag{A.47}$$

Analogously to (A.45), we also have for $m$ and $n$ large enough,

$$\begin{aligned}\sum_{i=n+1}^{m}\sum_{j=1}^{n}E[v_i^2v_j^2] &= \sum_{i=n+1}^{m}\sum_{j=1}^{n}[C_1(i,j)+C_2(i,j)]\\ &= \sigma_{m-n}^2\sigma_n^2(1+o(1)).\end{aligned} \tag{A.48}$$

Therefore, (A.37)–(A.48) imply that as $m,n\to\infty$

$$\begin{aligned}E[J_{mn}^2] &= E\left[\frac{1}{\sigma_m^2}\sum_{i=1}^{m}v_i^2-\frac{1}{\sigma_n^2}\sum_{j=1}^{n}v_j^2\right]^2\\ &= E\left[\frac{1}{\sigma_m^2}\sum_{i=n+1}^{m}v_i^2+\frac{\sigma_n^2-\sigma_m^2}{\sigma_m^2\sigma_n^2}\sum_{j=1}^{n}v_j^2\right]^2\\ &= \frac{1}{\sigma_m^4}E\left[\sum_{i=n+1}^{m}v_i^2\right]^2+\frac{(\sigma_n^2-\sigma_m^2)^2}{\sigma_m^4\sigma_n^4}E\left[\sum_{j=1}^{n}v_j^2\right]^2\\ &\quad -2\frac{\sigma_m^2-\sigma_n^2}{\sigma_m^4\sigma_n^2}\sum_{i=n+1}^{m}\sum_{j=1}^{n}E[v_i^2v_j^2]\\ &= \left(\frac{(m-n)^3}{m^3}+\frac{(m^{3/2}-n^{3/2})^2}{m^3}-2\frac{(m^{3/2}-n^{3/2})(m-n)^{3/2}}{m^3}\right)\\ &\quad \times(1+o(1))\\ &\to (1-r)^3+(1-r^{3/2})^2-2(1-r^{3/2})(1-r)^{3/2}\\ &= ((1-r)^{3/2}-(1-r^{3/2}))^2\ge 0\end{aligned}$$

using $\sigma_m^2=\frac{16J_{02}}{3\sqrt{2\pi}}m^{3/2}h$, $\sigma_n^2=\frac{16J_{02}}{3\sqrt{2\pi}}n^{3/2}h$ and $r=\lim_{m,n\to\infty}\frac{n}{m}$.

Since $r=0$ from the construction in the beginning of the proof of Lemma A.3 above, we have therefore shown the second part of (A.35). We now turn to the first part of (A.35). Using the results that $\sum_{i=n+1}^{m}E[v_i^4]=o(\sigma_{m-n}^4)$ and $\sum_{j=1}^{n}E[v_j^4]=o(\sigma_n^4)$, the proof of the first part of (A.35) follows from (A.36). We therefore have completed the proof of Lemma A.5. □

Define a random variable $N(T)$ in the same way as $T(n)$ that is defined in Karlsen and Tjøstheim (2001) [see Appendix B of Gao et al. (2008) for more details]. Recall

$$C_{10}=\frac{16J_{02}}{3\sqrt{2\pi}}\quad\text{and}\quad \sigma_T^2=C_{10}T^{3/2}h. \tag{A.49}$$

LEMMA A.6. *Let the conditions of Theorem 2.2 hold. Then as* $T \to \infty$

$$\frac{\widetilde{\sigma}_T^2}{\sigma_T^2} \to_D \xi^2 \tag{A.50}$$

*with* $\xi^2 = \frac{\sqrt{\pi}}{2} M_{1/2}(1)$, *where* $M_{1/2}(\cdot)$ *is a special case of the Mittag–Leffer process* $M_\beta(\cdot)$ *for* $\beta = \frac{1}{2}$ *as described by Karlsen and Tjøstheim (2001), page 388.*

PROOF. Observe that

$$\widetilde{\sigma}_T^2 = 2\sum_{t=1}^T \left(\sum_{s=1, s\neq t}^T a_{st}^2 u_s^2\right) u_t^2 = 2\sum_{t=1}^T \left(\sum_{s=1}^T a_{st}^2 u_s^2\right) u_t^2 - 2\sum_{t=1}^T a_{tt}^2 u_t^4.$$

Similarly to computations made between (A.5) and (A.6), it can be shown that

$$E\left[\sum_{t=1}^T \left(\sum_{s=1, s\neq t}^T a_{st}^2 (u_s^2 - 1)\right) u_t^2\right]^2 = o(\sigma_T^4) \tag{A.51}$$

using $E[u_1^2] = 1$.

Let $Q(u) = \frac{K^2(u)}{J_{02}}$. Then $Q(\cdot)$ is a probability kernel. Applying Lemma C.1 in Appendix C of Gao et al. (2008), we may show that as $T \to \infty$

$$\begin{aligned} &\frac{1}{T}\sum_{t=1}^T \left(\frac{1}{N(T)h}\sum_{s=1}^T Q\left(\frac{X_{s-1} - X_{t-1}}{h}\right)\right) u_t^2 \\ &\quad = \frac{1}{T}\sum_{t=1}^T \left(\frac{1}{N(T)h}\sum_{s=1}^T Q\left(\frac{X_{s-1} - X_{t-1}}{h}\right) - 1\right) u_t^2 \\ &\qquad + \frac{1}{T}\sum_{t=1}^T u_t^2 \to_P 1, \end{aligned} \tag{A.52}$$

where we have used the result that $\pi_s(Q) = \int Q(u)\,du \equiv 1$ [see the discussion at the end of Appendix B of Gao et al. (2008)].

Meanwhile, Theorem 3.2 of Karlsen and Tjøstheim (2001), page 389, is applicable to the current case of $X_t = X_{t-1} + u_t$ under $H_0$ to show that as $T \to \infty$

$$\frac{N(T)}{L_0\sqrt{T}} \to_D M_{1/2}(1), \tag{A.53}$$

when the slowly varying function $L_s(T)$ in this case is $L_s(T) \equiv L_0 = \frac{2\sqrt{2}}{3}$.

Thus, along with a strengthened version of Theorem 5.1 of Karlsen and Tjøstheim (2001), (A.51)–(A.53) imply as $T \to \infty$

$$\begin{aligned}
&\frac{2}{\sigma_T^2} \sum_{t=1}^{T} \left( \sum_{s=1}^{T} a_{st}^2 u_s^2 \right) u_t^2 \\
&\quad = \frac{2}{C_{10} T^{\frac{3}{2}} h} \sum_{t=1}^{T} \left( \sum_{s=1}^{T} a_{st}^2 u_s^2 \right) u_t^2 \\
&\quad = \frac{2L_0}{C_{10}} \frac{N(T)}{L_0 \sqrt{T}} \frac{1}{T} \sum_{t=1}^{T} \left( \frac{1}{N(T)h} \sum_{s=1}^{T} a_{st}^2 u_s^2 \right) u_t^2 \\
&\quad = \frac{2L_0}{C_{10}} \frac{N(T)}{L_0 \sqrt{T}} \frac{1}{T} \sum_{t=1}^{T} \left( \frac{1}{N(T)h} \sum_{s=1}^{T} a_{st}^2 (u_s^2 - 1) \right) u_t^2 \\
&\qquad + \frac{2L_0 J_{02}}{C_{10}} \frac{N(T)}{L_0 \sqrt{T}} \frac{1}{T} \sum_{t=1}^{T} \left( \frac{1}{N(T)h} \sum_{s=1}^{T} Q\left( \frac{X_{s-1} - X_{t-1}}{h} \right) \right) u_t^2 \\
&\quad \to_D \frac{\sqrt{\pi}}{2} M_{1/2}(1) \equiv \xi^2,
\end{aligned} \tag{A.54}$$

where we have used the facts that $\{u_s\}$ is a sequence of i.i.d. random errors with $E[u_1] = 0$ and $E[u_1^2] = 1$ and that $\{a_{st}^2 u_s^2 : 1 \leq s \leq t-1\}$ is independent of $u_t$. Therefore, (A.51)–(A.54) complete the proof of Lemma A.6. □

**Acknowledgments.** The authors would like to thank the coeditor, the Associate Editors and referees for their encouragements and constructive suggestions and comments. Thanks also go to several conference and seminar participants, in particular, Bruce Hansen, Yongmiao Hong, Qi Li, Joon Park, Peter Phillips, Robert Taylor and Qiying Wang, for offering their insightful suggestions and comments.

We would also like to thank Gowry Sriananthakumar and Jiying Yin for their excellent computing assistance.

## REFERENCES

Arapis, M. and Gao, J. (2006). Empirical comparisons in short-term interest rate models using nonparametric methods. *J. Financial Econometrics* **4** 310–345.

Bandi, F. and Phillips, P. C. B. (2003). Fully nonparametric estimation of scalar diffusion models. *Econometrica* **71** 241–283. MR1956859

Brockwell, P. and Davis, R. (1990). *Time Series Theory and Methods*. Springer, New York. MR0868859

Chow, Y. S. and Teicher, H. (1988). *Probability Theory*. Springer, New York. MR0953964

Dickey, D. A. and Fuller, W. A. (1979). Distribution of estimators for autoregressive time series with a unit root. *J. Amer. Statist. Assoc.* **74** 427–431. MR0548036

Fan, J. and Yao, Q. (2003). *Nonlinear Time Series: Nonparametric and Parametric Methods*. Springer, New York. MR1964455

Fan, Y. and Linton, O. (2003). Some higher theory for a consistent nonparametric model specification test. *J. Statist. Plann. Inference* **109** 125–154. MR1946644

Gao, J. (2007). *Nonlinear Time Series: Semiparametric and Nonparametric Methods*. Chapman & Hall, London. MR2297190

Gao, J. and Gijbels, I. (2008). Bandwidth selection in nonparametric kernel testing. *J. Amer. Statist. Assoc.* **484** 1584–1594.

Gao, J. and King, M. L. (2004). Adaptive testing in continuous-time diffusion models. *Econometric Theory* **20** 844–883. MR2089144

Gao, J., King, M. L., Lu, Z. and Tjøstheim, D. (2006). Specification testing in nonlinear time series with nonstationarity. Available at `http://www.adelaide.edu.au/directory/jiti.gao`.

Gao, J., King, M. L., Lu, Z. and Tjøstheim, D. (2008). Specification testing in nonlinear and nonstationary time series autoregression. Supplementary material available from the first author.

Gao, J., Lu, Z. and Tjøstheim, D. (2006). Estimation in semiparametric spatial regression. *Ann. Statist.* **34** 1395–1435. MR2278362

Granger, C. W. J., Inoue, T. and Morin, N. (1997). Nonlinear stochastic trends. *J. Econometrics* **81** 65–92. MR1484582

Granger, C. W. J. and Teräsvirta, T. (1993). *Modelling Nonlinear Dynamic Relationships*. Oxford Univ. Press, Oxford.

Hall, P. and Heyde, C. (1980). *Martingale Limit Theory and Its Applications*. Academic Press, New York. MR0624435

Hjellvik, V., Yao, Q. and Tjøstheim, D. (1998). Linearity testing using local polynomial approximation. *J. Statist. Plann. Inference* **68** 295–321. MR1629587

Jiang, G. (1998). Nonparametric modelling of US interest rate term structure dynamics and implication on the prices of derivative securities. *J. Financial and Quantitative Analysis* **33** 465–497.

Karlsen, H. and Tjøstheim, D. (1998). *Nonparametric Estimation in Null Recurrent Time Series* **50**. Sonderforschungsbereich 373, Humboldt Univ., Berlin.

Karlsen, H. and Tjøstheim, D. (2001). Nonparametric estimation in null recurrent time series. *Ann. Statist.* **29** 372–416. MR1863963

Karlsen, H., Myklebust, T. and Tjøstheim, D. (2007). Nonparametric estimation in a nonlinear cointegration model. *Ann. Statist.* **35** 252–299. MR2332276

Li, Q. (1999). Consistent model specification tests for time series econometric models. *J. Econometrics* **92** 101–147. MR1706996

Li, Q. and Wang, S. (1998). A simple consistent bootstrap tests for a parametric regression functional form. *J. Econometrics* **87** 145–165. MR1648892

Masry, E. and Tjøstheim, D. (1995). Nonparametric estimation and identification of nonlinear ARCH models. *Econometric Theory* **11** 258–289. MR1341250

Park, J. and Phillips, P. C. B. (2001). Nonlinear regressions with integrated time series. *Econometrica* **69** 117–162. MR1806536

Phillips, P. C. B. and Park, J. (1998). Nonstationary density estimation and kernel autoregression. Cowles Foundation Discussion Paper, No. 1181, Yale Univ.

Tong, H. (1990). *Nonlinear Time Series: A Dynamical System Approach*. Oxford Univ. Press, Oxford. MR1079320

Wang, Q. and Phillips, P. C. B. (2009). Asymptotic theory for local time density estimation and nonparametric cointegrating regression. *Econometric Theory* **25** 710–738.

J. Gao
School of Economics
University of Adelaide
Adelaide SA 5005
Australia
E-mail: jiti.gao@adelaide.edu.au

M. King
Office of Deputy Vice–Chancellor (Research)
Monash University
Melbourne VIC 3800
Australia
E-mail: max.king@adm.monash.edu.au

Z. Lu
Department of Mathematics
and Statistics
Curtin University of Technology
Perth WA 6450
Australia
E-mail: z.lu@curtin.edu.au

D. Tjøstheim
Department of Mathematics
University of Bergen
Bergen 5007
Norway
E-mail: dagt@mi.uib.no